\documentclass[12pt]{article}
\usepackage{amsfonts,amssymb,amsmath,fullpage,rotating}
\newcommand{\R}{{\mathbb R}}
\newcommand{\Z}{{\mathbb Z}}
\newcommand{\C}{{\mathbb C}}
\newcommand{\mH}{{\cal Q}}

\newcommand{\T}{{\mathbb T}}
\newcommand{\I}{{\mathbb I}}
\newcommand{\V}{{\mathbb V}}
\newcommand{\mO}{{\mathbb O}}
\newcommand{\D}{{\mathbb D}}

\newcommand{\mmod}{{\rm mod}\,}
\newcommand{\Fix}{{\rm Fix}\,}
\newcommand{\ri}{{\rm i}}
\newcommand{\rd}{{\rm d}}
\newcommand{\re}{{\rm e}}

\newcommand{\br}{\bf r}
\newcommand{\bz}{\bf z}

\newcommand{\bl}{\bf l}
\newcommand{\rl}{{\,|\,}}

\newtheorem{definition}{Definition}
\newtheorem{theorem}{Theorem}
\newtheorem{lemma}{Lemma}
\newtheorem{proposition}{Proposition}

\newtheorem{remark}{Remark}
\newcommand{\proof}{\noindent{\bf Proof: }}

\newcommand{\qed}{\hfill{\bf QED}\vspace{5mm}}

\begin{document}
\title{Asymptotic stability\\
of pseudo-simple heteroclinic cycles in $\R^4$}

\author{Olga Podvigina\\
Institute of Earthquake Prediction Theory\\
and Mathematical Geophysics\\
84/32 Profsoyuznaya St, 117997 Moscow, Russian Federation;\\
and\\
Pascal Chossat\\
Laboratoire J-A Dieudonn\'e, CNRS - UNS\\
Parc Valrose\\
06108 Nice Cedex 2, France
}

\maketitle

\begin{abstract}
Robust heteroclinic cycles in equivariant dynamical systems in $\R^4$
have been a subject of intense scientific investigation because, unlike
heteroclinic cycles in $\R^3$, they can have an intricate geometric structure
and complex asymptotic stability properties that are not yet completely
understood. In a recent work, we have compiled an exhaustive list of finite
subgroups of O(4) admitting the so-called {\em simple} heteroclinic cycles,
and have identified a new class which we have called {\em pseudo-simple}
heteroclinic cycles. By contrast with simple heteroclinic cycles,
a pseudo-simple one has at least one equilibrium with an unstable manifold
which has dimension 2 due to a symmetry. Here, we analyse the dynamics of nearby
trajectories and asymptotic stability of pseudo-simple
heteroclinic cycles in $\R^4$.
\end{abstract}

\section{Introduction}\label{sec_intro}

It is known since the 80's that vector fields defined on a vector space $V$
(or more generally a Riemannian manifold), which commute with the action of
a group $\Gamma$ of isometries of $V$, can possess invariant sets which are
structurally stable (within their $\Gamma$ symmetry class) and which are
composed of a sequence of saddle equilibria $\xi_1,\dots, \xi_M$ and
a sequence of trajectories $\kappa_i$, such that $\kappa_i$ belong to the unstable
manifold of $\xi_i$ as well as to the stable manifold of $\xi_{i+1}$ for all
$i=1,\dots,M$ where $\xi_{M+1}=\xi_1$. These objects are called {\em robust
heteroclinic cycles} and under certain conditions they can be asymptotically
attracting.
Robust heteroclinic cycles have been considerably studied in low-dimensional
vector spaces \cite{Kru97,cl2000}. Although their properties and classification
are well established in $V=\R^3$, the case of $\R^4$ appears to be much richer
and yet not completely investigated.
In \cite{km95a} the so-called {\it simple} robust heteroclinic
cycles were introduced, which can be defined as follows.
Let $\Delta_j$ denote the isotropy subgroup of an equilibrium $\xi_j$. The
heteroclinic cycle is simple if (i) the fixed point subspace of $\Delta_j$,
denoted $\Fix(\Delta_j)$, is an axis; (ii) for each $j$, $\xi_j$ is
a saddle and $\xi_{j+1}$ is a sink in an invariant plane $\Fix(\Sigma_j)$
(hence $\Sigma_j\subset\Delta_j\cap\Delta_{j+1}$); (iii) the isotypic
decomposition\footnote{The isotypic decomposition of the representation of a
group is the (unique) decomposition in a sum of equivalence classes of
irreducible representations.} for the action of $\Delta_j$ in the tangent space
at $\xi_j$ only contains one-dimensional components. In \cite{pc13} we have
found all finite groups $\Gamma\subset$\,O(4) {\em admitting} simple
heteroclinic cycles (i.e. for which there exists an open set of
$\Gamma$-equivariant vector fields possessing such an invariant set). We have
also pointed out that the definition of a simple heteroclinic cycle in former works
\cite{km95a,km04} had omitted condition (iii), which was implicitly assumed.
If this condition is not satisfied, then the behaviour of the
heteroclinic cycle turns out to be more complex, because there exists at least
one equilibrium $\xi_i$ at which the unstable manifold is two dimensional and,
moreover, is invariant under a faithful action of the dihedral group $\D_k$ for
some $k\geq 3$. We called heteroclinic cycles satisfying (i) and (ii) but not
(iii) {\em pseudo-simple}.

The aim of the present work is to investigate the dynamical properties of
pseudo-simple heteroclinic cycles in $\R^4$. We shall therefore consider
equations of the form
\begin{equation}\label{eq_ode}
\dot{\bf x}=f({\bf x}),\hbox{ where }f(\gamma{\bf x})=\gamma f({\bf x})
\hbox{ for all }\gamma\in\Gamma,\ \Gamma\subset\hbox{O(4)}\hbox{ is finite}
\end{equation}
and $f$ is a smooth map in $\R^4$, and which
possess a pseudo-simple cycle.

Our main result stated in Theorem 1, Section \ref{sec:instabSO(4)}, is that
if $\Gamma\subset$ SO(4), then the pseudo-simple cycle is {\em completely
unstable}. Namely, there exists a neighborhood of the heteroclinic cycle such
that any solution with initial condition in this neighborhood, except in a
subset of zero measure, leaves it in a finite time \cite{op12}.
Then we aim at analysing more deeply what the asymptotic dynamics can be in a
neighbourhood of the heteroclinic cycle.
We do this by focusing on specific examples. In Section \ref{sec4} we
introduce a subgroup of SO(4) which is algebraically elementary, however
possessing enough structure to allow for the existence of pseudo-simple
heteroclinic cycles. This group is a (reducible) representation of the
dihedral group $\D_3$. Hence, $\Gamma$ is algebraically isomorphic to $\D_3$.
In this example, the cycle involves two $\D_3$ invariant equilibria,
$\xi_1$ and $\xi_2$, such that for both the linearization $df(\xi_i)$ possess
a double eigenvalue with a $\D_3$ invariant associated eigenspace, one being negative and
the other being positive (unstable). We show that when the unstable double
eigenvalue is small, an attracting periodic orbit can exist in the vicinity of
the $\Gamma$-orbit of the heteroclinic cycle, the distance between the two
invariant sets vanishes as the unstable eigenvalue tends to zero.
This result is numerically illustrated by building an explicit third order
system for which a pseudo-simple cycle exists. In the next Section
\ref{sec:nonSO(4)} we extend the group $\Gamma$ by adding a generator which is
a reflection in $\R^4$. The resulting group is isomorphic to $\D_3\times \Z_2$
and possesses the same properties for the existence of a pseudo-simple cycle.
We show that in this case the cycle is not completely unstable anymore,
but {\em asymptotically fragmentarily stable} \cite{op12} instead. This means
that the cycle is not asymptotically stable, however
in its any small neighbourhood there exists a set of positive measure,
such that any solution with initial condition in that set remains in this
neighbourhood for all $t>0$ and converges to the cycle as $t\to\infty$.
This result is then illustrated numerically. \\
In Section \ref{sec:otherexample} we present another example
of a group in SO(4) admitting pseudo-simple cycles with
a more complex group structure and acting irreducibly on $\R^4$.
The Theorem of existence of nearby stable periodic orbits proven
in Section \ref{sec4} holds true for the same reasons as in the $\D_3$ case.
We then investigate numerically the asymptotic behavior of trajectories
near the cycle and show that periodic orbits of various kinds can be observed
as well as more complex dynamics when the unstable double eigenvalue becomes
larger. \\
In the last Section of this paper we discuss our findings and possible
continuation of the study.

\section{Pseudo-simple heteroclinic cycles} \label{sec:pseudosimple}
In this Section we introduce the notion of a pseudo-simple heteroclinic cycle and notations which will be used throughout the paper.

\subsection{Basic definitions}\label{sec:definitions}
Given a finite group $\Gamma$ of isometries of $\R^4$, let $\Sigma$ be an
{\em isotropy subgroup} of $\Gamma$ (the subgroup of elements in $\Gamma$ which
fix a point in $\R^4$). We denote by $\Fix(\Sigma)$ the subspace comprised of all
points in $\R^4$, which are fixed by $\Sigma$. \\
Note that (i) if $\Sigma\subset\Delta\subset\Gamma$, then $\Fix(\Delta)\subset \Fix(\Sigma)$; (ii) $\gamma\cdot \Fix(\Sigma)=\Fix(\gamma\Sigma\gamma^{-1})$ for any $\gamma\in\Gamma$. \\
In the following we shall assume the $\Gamma$-equivariant system
(\ref{eq_ode}) admits a sequence of isotropy subgroups $\Delta_j\subset\Gamma$,
$j=1,\dots,M$, such that the following holds:
\begin{itemize}
\item[(i)] dim$\Fix(\Delta_j)=1$ (axis of symmetry). We denote $L_j=\Fix(\Delta_j)$.
\item[(ii)] On each $L_j$, there exists an equilibrium
$\xi_j$ of \eqref{eq_ode}.
\item[(ii)] For each $j$, $L_j=P_{j-1}\cap P_j$, where $P_j$ is a plane of
symmetry: $P_j=\Fix(\Sigma_j)$ for some $\Sigma_j\subset\Delta_j$. We set
$P_M=P_0$.
\item[(iv)] $\xi_j$ is a sink in $P_{j-1}$ and a saddle in $P_j$. Moreover,
a saddle-sink heteroclinic trajectory $\kappa_j$ of \eqref{eq_ode} connects
$\xi_j$ to $\xi_{j+1}$ in $P_j$ ($j=1,\dots,M$). Note that $\kappa_M$ connects
$\xi_M$ to $\xi_1$ in $P_M$.
\end{itemize}
Conditions (i) - (iv) insure the existence of a robust heteroclinic cycle for \eqref{eq_ode}.
We denote by $J_j$ the Jacobian matrix $df(\xi_j)$. Since $L_j$ is
flow-invariant, $J_j$ has a {\em radial} eigenvalue $r_j$ with eigenvector along
$L_j$ and by (iv) we have $r_j<0$. Moreover (iv) implies that in $P_{j-1}$, $J_j$
has a {\em contracting} eigenvalue $-c_j$ ($c_j>0$), while in $P_j$ it has an
{\em expanding} eigenvalue $e_j>0$ corresponding to the unstable eigendirection
in that plane. The remaining eigenvalue $t_j$ is called {\em transverse},
its associated eigenspace is the complement to $P_{j-1}\oplus P_j$ in $\R^4$.

\begin{definition} \label{def:admits}
We say that the group $\Gamma$ {\em admits} robust heteroclinic
cycles if there exists an open subset of the set of smooth $\Gamma$-equivariant
vector fields in $\R^n$, such that vector fields in this subset possess a
(robust) heteroclinic cycle.
\end{definition}

In order to insure the existence of a robust heteroclinic cycle
it is enough to find $m\leq M$ and $\gamma\in\Gamma$ such that a minimal sequence of robust
heteroclinic connections $\xi_1\rightarrow\cdots\rightarrow\xi_{m+1}$ exists with $\xi_{m+1}=\gamma\xi_1$
(minimal in the sense that no equilibria in this sequence satisfy $\xi_i=\gamma'\xi_j$
for some $\gamma'\in\Gamma$ and $1\le i,j\le m$). It follows that $\gamma^k=1$ where $k$ is a divisor of $M$.
\begin{definition}\label{def42}
The sequence $\xi_1\rightarrow\cdots\rightarrow\xi_{m}$ and the element $\gamma$ define
a {\em building block} of the heteroclinic cycle.
\end{definition}

The asymptotic dynamics in a neighborhood of a robust heteroclinic cycle is what
makes these objects special. A heteroclinic cycle is called {\em asymptotically
stable} if it attracts all trajectories in its small neighborhood.
When a heteroclinic cycle is not asymptotically stable, it can however possess
residual types of stability. Below we define the notions, which are relevant
in the context of this paper, see \cite{op12}. Given a heteroclinic cycle $X$ and writing the flow of \eqref{eq_ode} as $\Phi_t(x)$, the $\delta$-basin of attraction of $X$ is the set
$$
B_\delta(X) = \{x\in \R^4~;~d(\Phi_t(x),X)<\delta\hbox{~for all~}t>0\mbox{~and~} \lim_{t\rightarrow+\infty}d(\Phi_t(x),X)=0 \}.
$$
\begin{definition} \label{def:completelyunstable}
A heteroclinic cycle $X$ is {\em completely unstable} if there exists $\delta>0$
such that $B_\delta(X)$ has Lebesgue measure 0.
\end{definition}
\begin{definition} \label{def:fragmstable}
A heteroclinic cycle $X$ is {\em fragmentarily asymptotically stable} if for any $\delta>0$, $B_\delta(X)$ has positive Lebesgue measure.
\end{definition}

\subsection{Simple and pseudo-simple heteroclinic cycles}
Recall that the isotypic decomposition of a representation $T$ of a (finite)
group $G$ in a vector space $V$ is the decomposition
$V=V^{(1)}\oplus\cdots\oplus V^{(r)}$ where $r$ is the number of equivalence
classes of irreducible representations of $G$ in $V$ and each
$V^{(j)}=T_{|V_j}$ is the sum of the equivalent irreducible representations in
the $j$-th class. This decomposition is unique. The subspaces $V^{(j)}$ are
mutually orthogonal (if $G$ acts orthogonally). \\
Then the following holds \cite{pc13}.

\begin{lemma}\label{prop:isotypic decomp}
Let a robust heteroclinic cycle in $\R^4$ be such that for all $j$:
(i) $\dim P_j=2$, (ii) each connected component of $L_j\setminus\{0\}$ is
intersected at most at one point by the heteroclinic cycle.
Then the isotypic decomposition of the
representation of $\Delta_j$ in $\R^4$ is of one of the following types:
\begin{enumerate}
\item $L_j\oplus^\perp V_j\oplus^\perp W_j\oplus^\perp T_j$ (the symbol $\oplus^\perp$ indicates the orthogonal direct sum).
\item $L_j\oplus^\perp V_j\oplus^\perp \widetilde{W_j}$ where $\widetilde{W_j}=W_j\oplus T_j$ has dimension 2.
\item $L_j\oplus^\perp \widetilde{V_j}\oplus^\perp W_j$ where $\widetilde{V_j}=V_j\oplus T_j$ has dimension 2.
\end{enumerate}
In cases 2 and 3, $\Delta_j$ acts in $\widetilde{W_j}$ (respectively, $\widetilde{V_j}$) as a dihedral group $\D_m$ in $\R^2$ for some $m\geq 3$. It
follows that in case 2, $e_j$ is double (and $e_j=t_j$) while in case 3, $-c_j$ is double (and $-c_j=t_j$).
\end{lemma}

\begin{definition}\label{def:verysimple}
A robust heteroclinic cycle in $\R^4$ satisfying the conditions 1 and 2
of Lemma \ref{prop:isotypic decomp} is called {\em simple} if case 1 holds true for all $j$, and
{\em pseudo-simple} otherwise.
\end{definition}

\begin{remark}\label{rem:Zk}
Definition \ref{def:verysimple} implies that any pseudo-simple heteroclinic
cycle in a $\Gamma$-equivariant system (\ref{eq_ode}),
where $\Gamma\subset$\,SO(4), has at least one $1\le j\le M$ such that
$\Sigma_j\cong\Z_k$ with $k\ge3$.
\end{remark}

\begin{remark} \label{rem:planereflection}
An order two element $\sigma$ in SO(4) whose fixed point subspace is a plane $P$
must act as $-Id$ in the plane $P^\perp$ fully perpendicular to $P$. Nevertheless,
to distinguish it from other rotations fixing the points on $P$, we call $\sigma$ a \rm{plane reflection}.
\end{remark}

\section{Instability of pseudo-simple cycles with $\Gamma\subset SO(4)$}\label{sec:instabSO(4)}

In this Section we prove the following Theorem:
\begin{theorem}\label{th1}
Let $\xi_1\to\ldots\to\xi_M$ be a pseudo-simple heteroclinic cycle in
a $\Gamma$-equivariant system, where $\Gamma\subset$\,SO(4). Then generically
the cycle is completely unstable.
\end{theorem}
\proof
Following \cite{km95a,km04,op12,op13}, to study asymptotic stability of
a heteroclinic cycle, we approximate
a "first return map" on a transverse (Poincar\'e) section of the cycle and
consider its iterates for trajectories that stay in a small neighbourhood
of the cycle.\\
In Section \ref{sec:definitions} we have defined radial, contracting, expanding and transverse eigenvalues of the linearisation $J_j=df(\xi_j)$.
Let $(\tilde u,\tilde v,\tilde w,\tilde q)$
be coordinates in the coordinate system with the origin at $\xi_j$ and the basis
comprised of the associated eigenvectors in the following order: radial,
contracting, expanding and transverse.
If $\tilde\delta$ is small, in a $\tilde\delta$-neighbourhood of $\xi_j$ the system (\ref{eq_ode}) can be approximated by the linear system
\begin{equation}\label{lmap}
\begin{array}{l}
\dot u=-r_ju\\
\dot v=-c_j v\\
\dot w=e_j w\\
\dot q=t_jq.
\end{array}
\end{equation}
Here, $(u,v,w,q)$ denote the scaled coordinates
$(u,v,w,q)= (\tilde u,\tilde v,\tilde w,\tilde q)/\tilde\delta$. In fact a
version of the Hartman-Grobman theorem exists \cite{cl2000}, which allows to
linearize the system with a $C^k$ ($k\geq 1$) equivariant change of variables in a neighbourhod of $\xi_j$ if conditions of nonresonance are satisfied between the eigenvalues (which is a generic condition here).

Let $(u_0,v_0)$ be the point in $P_{j-1}$ where the trajectory $\kappa_{j-1}$
intersects with the circle $u^2+v^2=1$, and $h$ be local coordinate in the line
tangent to the circle at the point $(u_0,v_0)$. We consider two crossections:
$$
\widetilde H^{(in)}_j=\{(h,w,q)~:~|h|,|w|,|q|\le1\}
$$
and
$$
\widetilde H^{(out)}_j=\{(u,v,w,q)~:~|u|,|v|,|q|\le1, w=1\}.
$$

Near $\xi_j$, trajectories of (\ref{eq_ode}) with initial condition in
$\widetilde H^{(in)}_j$ hit the outgoing section $\widetilde H^{(out)}_j$
after a time $\tau$, which tends to infinity as the initial condition comes
closer to $P_{j-1}$.
The global map $\tilde\psi_j:\widetilde H^{(out)}_j\to\widetilde H^{(in)}_{j+1}$
associates a point where a trajectory crosses $\widetilde H^{(in)}_{j+1}$ with
the point where it previously crossed $\widetilde H^{(out)}_j$. It is defined
in a neighborhood of $0$ (in the local coordinates in $\widetilde H^{(out)}_j$)
and is a diffeomorphism. The first return map is now defined as the composition $\widetilde g=g_M\circ\ldots\circ g_1$, where $g_j=\tilde\psi_j\circ \tilde\phi_j$.

As it is shown in \cite{pa11,op12}, in the study of stability only the
coordinates that are in $P_j^{\perp}$ are of importance. By $\phi_j$ and
$\psi_j$ we denote the maps $\tilde\phi_j$ and $\tilde\psi_j$ restricted to
$P_j^{\perp}$ (to be more precise
$H^{(in)}_j=P_{j-1}^{\perp}\cap\widetilde H^{(in)}_j$ and
$H^{(out)}_j=P_j^{\perp}\cap\widetilde H^{(out)}_j$). \\
Now, by definition a pseudo-simple cycle has at least one connection $\kappa_j:\xi_j\to\xi_{j+1}$
such that $\Sigma_j\cong\Z_k$ with $k\ge3$ (see remark \ref{rem:Zk}).
We can assume that $j=1$. The equilibria
$\xi_1$ and $\xi_2$ belong to the plane $P_1=\Fix(\Sigma_1)$.
We now prove existence of $\varepsilon>0$ such that
\begin{equation}\label{nonincl}
\phi_2\circ\psi_1\circ\phi_1(H^{(in)}_1(\varepsilon))\cap H^{(out)}_2(\varepsilon)=\varnothing,
\end{equation}
where
$$H^{(in)}_1(\varepsilon)=\{(w,q)\in H^{(in)}_1~:~|(w,q)|<\varepsilon\}
\hbox{ and }H^{(out)}_2(\varepsilon)=
\{(v,q)\in H^{(out)}_2~:~|(v,q)|<\varepsilon\}.$$
Evidently, because of (\ref{nonincl}) the first return map in
$H^{(in)}_1(\varepsilon)$ cannot be completed in general, which implies that
the cycle is completely unstable. \\
To prove \eqref{nonincl} we employ polar coordinates
$(\rho_1,\theta_1)$ and $(\rho_2,\theta_2)$, in $H^{(out)}_1$ and $H^{(in)}_2$,
respectively, such that $v_1=\rho_1\cos\theta_1$, $q_1=\rho_1\sin\theta_1$,
$w_2=\rho_2\cos\theta_2$ and $q_2=\rho_2\sin\theta_2$.
Due to (\ref{lmap}) in the leading order the maps $\phi_1$ and $\phi_2$ are
\begin{equation}\label{phmap1}
(\rho_1,\theta_1)=\phi_1(w_1,q_1)=(v_{0,1}w_1^{c_1/e_1},\arctan(q_1/v_{0,1})),
\end{equation}
\begin{equation}\label{phmap2}
(v_2,q_2)=\phi_2(\rho_2,\theta_2)=
(v_{0,2}(\rho_2\cos\theta_2)^{c_2/e_2},\tan\theta_2).
\end{equation}

In the leading order the global map
$\psi_1:H^{(out)}_1\to H^{(in)}_2$ is linear.
The map commutes with the group $\Sigma_2$ acting on $P_2^{\perp}$ as rotation
by $2\pi/k$, where $k\ge3$, therefore in the leading order $\psi_1$ is a rotation by
a finite angle composed with a linear transformation of $\rho$,
\begin{equation}\label{glmap1}
(\rho_2,\theta_2)=\psi_1(\rho_1,\theta_1)=(A\rho_1,\theta_1+\Theta),
\end{equation}
where generically $\Theta\ne N\pi/k$.

Denote
\begin{equation}\label{setde0}
\alpha=\min_{1\ge N\ge 2k} |\Theta- N\pi/k|
\end{equation}
and set
\begin{equation}\label{setde}
0<\varepsilon<\min\biggl(\tan{\alpha\over2},~v_{0,1}\tan{\alpha\over2}\biggr).
\end{equation}
Any $(w_1,q_1)\in H^{(in)}_1(\varepsilon)$ satisfies $q_1<\varepsilon$,
therefore (\ref{phmap1}) and (\ref{setde}) imply that $\theta_1<\alpha/2$.
Hence, due to (\ref{glmap1})-(\ref{setde}) $|\theta_2- N\pi/k|>\alpha/2$ for any $k$.
The steady state $\xi_2$ has $k$ symmetric copies (under the action of
symmetries $\sigma\in\Sigma_2$) of the heteroclinic connection
$\kappa_2:\xi_2\to\xi_3$ which belong to the hyperplanes
$\theta_2=N\pi/k$ with some integer $N$'s. Due to (\ref{phmap2}) and (\ref{setde}),
the distance of $(v_2,q_2)$ to any of these hyperplanes is larger than
$\tan(\alpha/2)$, which implies (\ref{nonincl}).
\qed

\section{Existence of nearby periodic orbits when $\Gamma\subset SO(4)$:
an example}\label{sec4}

 We have proved that any pseudo-simple heteroclinic cycle
in a $\Gamma$-equivariant system, where $\Gamma\subset$\,SO(4), is completely unstable.
In this Section we show that nevertheless trajectories staying
in a small neighbourhood of a pseudo-simple cycle for all $t>0$ can exist.
Namely, we give an example of a subgroup of SO(4) admitting pseudo-simple
heteroclinic cycles and prove generic existence of an asymptotically stable
periodic orbit near a heteroclinic cycle if an unstable double eigenvalue
is sufficiently small.

\subsection{Definition of $\Gamma$ and existence of a pseudo-simple cycle}\label{exa1}

We write $(x_1,y_1,x_2,y_2)\in\R^4$ and $z_j=x_j+iy_j$. Let $\Gamma$ be the group generated by the transformations
\begin{equation} \label{gen:Gamma}
\rho:~(z_1,z_2)\mapsto (z_1,e^{\frac{2\pi i}{3}}z_2),~~\kappa(z_1,z_2)\mapsto (\bar z_1,\bar z_2)
\end{equation}
This group action obviously decomposes into the direct sum of three irreducible
representations of the dihedral group $\D_3$: \\
(i) the trivial representation acting on the component $x_1$, \\
(ii) the one-dimensional representation acting on $y_1$ by $\kappa y_1=-y_1$, \\
(iii) the two-dimensional natural representation of $\D_3$ acting on $z_2=(x_2,y_2)$.

There are two types of fixed-point subspaces and one type of invariant axis for this action: \\
(i) $P_1=\{(x_1,y_1,0,0)\}=\Fix(\rho)$, \\
(ii) $P_2=\{(x_1,0,x_2,0)\}=\Fix(\kappa)$, \\
(iii) $L=P_1\cap P_2=\Fix(\D_3)$. \\
Note that $P_1$ is $\Gamma$-invariant, while $P_2$ has two symmetric
copies, $P'_2=\rho P_2$ and $P''_2=\rho^2P_2$.

\begin{proposition}\label{prop:admissible}
The group $\Gamma$, generated by $\rho$ and $\kappa$ (\ref{gen:Gamma}),
{\em admits} pseudo-simple cycles (see Definition \ref{def:admits}).
\end{proposition}
\proof
The proof of existence of an open set of smooth $G$ equivariant vector fields
with saddle-sink orbits in $P_1$ and in $P_2$ connecting equilibria $\xi_1$ and $\xi_2$ lying on the axis $L$ and such that $\xi_1$ is a sink in $P_2$ and $\xi_2$ is a sink in $P_1$, goes along the same arguments as in Lemma 5 of \cite{pc13}. We just need to check that the cycle is pseudo-simple. This comes from the fact that the Jacobian matrix of the vector field taken at $\xi_1$ or $\xi_2$ has a double eigenvalue with eigenspace $\{(0,0,x_2,y_2)\}$ due to the action of $\rho$, which generates the subgroup $\Z_3$ (negative eigenvalue at $\xi_1$ and positive eigenvalue at $\xi_2$).
\qed

By letting $\Gamma$ act on the heteroclinic cycle, one obtains a 6-element orbit
of cycles which all have one heteroclinic connection in $P_1$ and
the other one in $P_2$, $P'_2$ or $P''_2$.

\subsection{Existence and stability of periodic orbits}
We prove the following Theorem.

\begin{theorem}\label{thperorb}
 Consider the $\Gamma$-equivariant system
\begin{equation}\label{eqal_ode}
\dot{\bf x}=f({\bf x},\mu),\hbox{ where }
f(\gamma{\bf x},\mu)=\gamma f({\bf x},\mu)\mbox{ for all }\gamma\in\Gamma,
\end{equation}
$f:\R^4\times\R\to\R^4$ is a smooth map and
$\Gamma$ is generated by $\rho$ and $\kappa$ (\ref{gen:Gamma}).
Let $\xi_1$ and $\xi_2$ be the equilibria introduced in proposition
\ref{prop:admissible},
$-c_j$ and $e_j$ be the eigenvalues of $df(\xi_j)$, such that $-c_1$ and $e_2$
are double eigenvalues with a natural action of $\D_3$ in their eigenspaces
and $e_1$ and $-c_2$ are the eigenvalues associated with the one-dimensional
eigenspace, where the action of $\Gamma$ in not trivial.
Suppose that there exists $\mu_0>0$ such that
\begin{itemize}
\item[(i)] $e_2<0$ for $-\mu_0<\mu<0$ and
$e_2>0$ for $0<\mu<\mu_0$;
\item[(ii)] for any $0<\mu<\mu_0$ there exist heteroclinic
connections $\kappa_1=(W_u(\xi_1)\cap P_1)\cap W_s(\xi_2)\ne\varnothing$ and
$\kappa_2=(W_u(\xi_2)\cap P_2)\cap W_s(\xi_1)\ne\varnothing$, introduced
in proposition \ref{prop:admissible}.
\end{itemize}
Denote by $X$ the group orbit of heteroclinic connections $\kappa_1$ and
$\kappa_2$:
$$X=(\cup_{\gamma\in\Gamma}\gamma\kappa_1)\bigcup
(\cup_{\gamma\in\Gamma}\gamma\kappa_2).$$
Then
\item[(a)] if $3c_1<e_1$ then there exist $\mu'>0$ and $\delta>0$, such that
for any $0<\mu<\mu'$ almost all trajectories escape from
$B_{\delta}(X)$ as $t\to\infty$;
\item[(b)] if $3c_1>e_1$ then generically there exists periodic orbit
bifurcating from $X$ at $\mu=0$.
To be more precise, for any $\delta>0$ we can find $\mu(\delta)>0$
such that for all $0<\mu<\mu(\delta)$ the system
(\ref{eqal_ode}) possesses an asymptotically stable periodic
orbit that belongs to $B_{\delta}(X)$.
\end{theorem}

The rest of this Section is devoted to the proof of the Theorem. We start with
two Lemmas which describe some properties of trajectories of a generic
$\D_3$-equivariant systems.

\subsubsection{Two Lemmas}

A generic $\D_3$-equivariant second order dynamical system in $\C$ is
\begin{equation}\label{sysA}
\dot z=\alpha z+\beta\bar z^2.
\end{equation}
In polar coordinates, $z=r\re^{\ri\theta}$, it takes the form
\begin{equation}\label{pcor1}
\begin{array}{rcl}
\dot r&=&\alpha r+\beta r^2\cos3\theta,\\
\dot\theta&=&-\beta r\sin3\theta.
\end{array}
\end{equation}
We assume that $\alpha>0$ and $\beta>0$.
The system has three invariant axes with $\theta=K\pi/3$, $K=0,1,2$, and three equilibria off
origin, with $r=\alpha/\beta$ and $\theta=(2k+1)\pi/3$, $k=0,1,2$.
We consider the system in the sector $0\le\theta<\pi/3$,
the complement part of $\C$ is related to this sector by symmetries of the group $\D_3$.

Trajectories of the system satisfy
\begin{equation}\label{pcor2}
{\rd r\over\rd\theta}=-{\alpha+\beta r\cos3\theta\over \beta\sin3\theta}.
\end{equation}
Re-writing this equation as
$${\rd r\over\rd\theta}+r{\cos3\theta\over\sin3\theta}=-{\alpha\over\beta\sin3\theta},$$
multiplying it by $\sin^{1/3}3\theta$ and integrating, we obtain that
\begin{equation}\label{trajr}
r\sin^{1/3}3\theta=-{\alpha\over\beta}S(\theta)+C,\hbox{ where }
S(\theta)=\int_0^{\theta}\sin^{-2/3}3\theta\rd\theta,
\end{equation}
which implies that
\begin{equation}\label{trajr1}
r\sin^{1/3}3\theta+{\alpha\over\beta}S(\theta)=
r_0\sin^{1/3}3\theta_0+{\alpha\over\beta}S(\theta_0)
\end{equation}
for the trajectory through the point $(r_0,\theta_0)$.
Let $\theta(r;r_0,\theta_0)$ be a solution to (\ref{pcor1}) with the
initial condition $r(0)=r_0$ and $\theta(0)=\theta_0$.
Then
$$C={\alpha\over\beta}S(\tilde\theta),\hbox{ where }
\tilde\theta=\lim_{r\to 0}\theta(r;r_0,\theta_0).$$
We can re-write (\ref{trajr}) as
\begin{equation}\label{trajro}
r\sin^{1/3}3\theta=-{\alpha\over\beta}S(\theta)+{\alpha\over\beta}S(\tilde\theta).
\end{equation}
Note, that
\begin{equation}\label{gamf}
S({\pi\over3})={\sqrt{\pi}\over 3}{\Gamma(1/6)\over\Gamma(2/3)},\ \hbox{ where }\Gamma
\hbox{ is the Gamma function.}
\end{equation}

\begin{lemma}\label{d3_1} Given $\beta>0$ and $\varepsilon>0$, there exists
$\alpha_0>0$ such that for any $0<\alpha<\alpha_0$, $0<r_0<\varepsilon$
and $0<\theta_0<\pi/3$
\begin{equation}\label{conda1}
r\sin\theta(r;r_0,\theta_0)<\varepsilon\hbox{ for all }t>0.
\end{equation}
\end{lemma}

\begin{figure}

\centerline{\includegraphics[width=12cm]{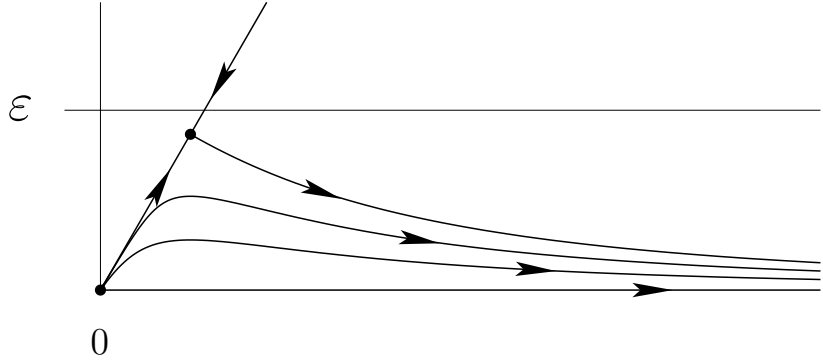}}

\vspace*{-42mm}
\hspace*{30mm}{\LARGE$\varepsilon$}

\vspace*{26mm}
\hspace*{41mm}{\large 0}

\vspace*{4mm}
\noindent
\caption{Trajectories of the system (\ref{pcor1}) in the sector $0\le\theta\le\pi/3$.}
\label{fig1}\end{figure}

\proof
Set $\alpha_0=\varepsilon\beta/4$. From (\ref{pcor1})
\begin{equation}\label{lem1_1}
{\rd\over \rd t}r\sin\theta=\dot r\sin\theta+\dot\theta r\cos\theta=
r\sin\theta(\alpha-\beta r\cos\theta/2).
\end{equation}
From (\ref{trajr}), for $0<\theta_0<\pi/3$ we have
$\lim_{r\to\infty}\theta(r;r_0,\theta_0)=0$, which implies that
$\rd/\rd t(r\sin\theta)<0$ for large $r$. Therefore, the maximum of
$r\sin\theta$ is achieved either at $t=0$, or at the point where
$2\alpha=\beta r\cos\theta$.
Since $0<\theta<\pi/3$, at this point
$r\sin\theta=2\alpha\sin\theta/(\beta\cos\theta)<4\alpha/\beta<\varepsilon$.
Since $r_0<\varepsilon$, we have $r\sin\theta<\varepsilon$ at $t=0$ as well.
(The behaviour of trajectories of the system (\ref{pcor1}) is shown in
Fig.~\ref{fig1}.)
\qed

\begin{lemma}\label{d3_30}
Let $\tau(r_0,\theta_0)$ denotes the time it takes the trajectory of the
system (\ref{pcor1}) starting at $(r_0,\theta_0)$ to reach $r=1$ and
$\vartheta(r_0,\theta_0)$ denotes the value of $\theta$ at $r=1$. Then
\begin{itemize}
\item[(i)] $\tau(r_0,0)$ satisfies
$${\rm e}^{\alpha\tau(r_0,0)}={r_0+\alpha/\beta\over r_0(1+\alpha/\beta)}.$$
\item[(ii)] $\tau(r_0,\theta_0)$ satisfies
\begin{equation}\label{estr}
\tau(r_0,\theta_0)>\tau(r_0,0)\hbox{ for any }0<\theta_0<\pi/3.
\end{equation}
\item[(iii)] $\vartheta(r_0,\theta_0)$ satisfies
\begin{equation}\label{estt}
\sin^{1/3}3\vartheta(r_0,\theta_0)+
{\alpha\over\beta}S(\vartheta(r_0,\theta_0))=
r_0\sin^{1/3}3\theta_0+{\alpha\over\beta}S(\theta_0).
\end{equation}
\item[(iv)] Given $C>0$, $\beta>0$ and $0<\theta_0<\pi/3$, for sufficiently
small $\alpha$ and $r_0$
$${\rm e}^{-C\tau(r_0,\theta_0)}\ll\vartheta(r_0,\theta_0).$$
\item[(v)]
For small $\alpha$, $r_0\ll\alpha$, $r'\ll r_0$ and $\theta'\ll1$
$$|\re^{-C\tau(r_0+r',\theta_0+\theta')}-
\re^{-C\tau(r_0,\theta_0)}|<\re^{-C\tau(r_0,\theta_0)}|
{r'\over\alpha r_0}+3\theta' S(\theta_0)\tau(r_0,\theta_0)|.$$
\end{itemize}
\end{lemma}

\proof
To obtain (i), we integrate the first equation in (\ref{pcor1}), where we set
$\theta=0$. Since in (\ref{pcor1}) $\theta<\theta'$ implies
$\dot r(r,\theta)>\dot r(r,\theta')$, (ii) holds true.
The equality (iii) follows from (\ref{trajr1}).

\medskip
Since $\dot\theta<0$ (see (\ref{pcor1})\,), (iii) implies that
$\vartheta(r_0,\theta_0)\ll1$ for small $\alpha$ and $r_0$. Therefore
\begin{equation}\label{esthet}
\vartheta(r_0,\theta_0)\approx{1\over3}
\biggl(r_0\sin^{1/3}3\theta_0+{\alpha\over\beta}S(\theta_0)\biggr)^3>
{1\over 3\beta^3}\biggl(\max(\alpha,\beta r_0)\min
(\sin^{1/3}3\theta_0,S(\theta_0))\biggr)^3.
\end{equation}
Due to (ii), to prove (iv) it is sufficient to show that
$${\rm e}^{-C\tau(r_0,0)}\ll\vartheta(r_0,\theta_0).$$
Because of (i), for small $r_0$ and $\alpha$
$${\rm e}^{-C\tau(r_0,0)}=
\biggl(
{r_0(1+\alpha/\beta)\over r_0+\alpha/\beta}\biggr)^{C/\alpha}\approx
[(1-s)^{1/s}]^{C/(\alpha+\beta r_0)},~\hbox{ where }~
s={\alpha\over\alpha+\beta r_0}.$$
Since $\sup_{0<s<1}(1-s)^{1/s}=\re^{-1}$, we have
\begin{equation}\label{estau}
{\rm e}^{-C\tau(r_0,0)}<\re^{-C/\max(\alpha,\beta r_0)}.
\end{equation}
Together, (\ref{esthet}) and (\ref{estau}) imply (iv).

\medskip
For $r\ll\alpha$ the first equation in (\ref{pcor1}) can be approximated by
$\dot r=\alpha r$. Therefore, the difference $\Delta\tau(r')=
\tau(r_0+r',\theta_0)-\tau(r_0,\theta_0)$
satisfies
$$r_0\re^{\alpha\Delta\tau(r')}\approx(r_0+r'),$$
which implies $\Delta\tau(r')\approx r'/\alpha r_0$. Hence,
\begin{equation}\label{tau1}
\re^{-C\tau(r_0+r',\theta_0)}-\re^{-C\tau(r_0,\theta_0)}\approx
\re^{-C\tau(r_0,\theta_0)}{r'\over\alpha r_0}.
\end{equation}

Consider two nearby trajectories, $\theta(r;r_0,\theta_0)$ and
$\theta(r;r_0,\theta_0+\theta')$. From (\ref{pcor1}),\break
$\theta(r;r_0,\theta_0+\theta')-\theta(r;r_0,\theta_0)<\theta'$ for all $t>0$,
therefore in view of (iii)
$${|\dot r(r,\theta(r;r_0,\theta_0+\theta'))-\dot r(r,\theta(r;r_0,\theta_0))|
\over |\dot r(r,\theta(r;r_0,\theta_0))|}
<{|3\theta'\beta r^2\sin3\theta|
\over |\alpha r+\beta r^2\cos3\theta|}<
|{3\beta\over\alpha}\theta'r\sin^{1/3}3\theta|<|3\theta' S(\theta_0)|,$$
which implies that
$$|\tau(r_0,\theta_0+\theta')-\tau(r_0,\theta_0)|<|3\theta' S(\theta_0)\tau(r_0,\theta_0)|.$$
Hence,
\begin{equation}\label{tau2}
|\re^{-C(\tau(r_0,\theta_0+\theta'))}-\re^{-C\tau(r_0,\theta_0)}|
<|\re^{-C\tau(r_0,\theta_0)}3\theta' S(\theta_0)\tau(r_0,\theta_0)|.
\end{equation}
Together (\ref{tau1}) and (\ref{tau2}) imply (v).
\qed

\subsubsection{Proof of Theorem \ref{thperorb}}
\medskip

As in the proof of Theorem \ref{th1}, we approximate trajectories in the
vicinity of the cycle by superposition of local and global maps.
We consider $g=\phi_1\psi_2\phi_2\psi_1:\ H^{(out)}_1\to H^{(out)}_1$,
where the map $\phi_1$ is given by (\ref{phmap1}).
The map $\psi_2$ in the leading order is linear:
\begin{equation}\label{psi2}
(w_1,z_1)=\psi_2(v_2,z_2)=(B_{11}v_2+B_{12}z_2,B_{21}v_2+B_{22}z_2).
\end{equation}

Because of (i), for small $\mu$ the expanding eigenvalue of $\xi_2$ depends
linearly on $\mu$, therefore without restriction of generality we can assume that
$e_2=\mu$.
Generically, all other eigenvalues and coefficients in the expressions for local
and global maps do not vanish for sufficiently small $\mu$ and are of the order
of one. We assume them to be constants independent of $\mu$. From (ii),
the eigenvalues satisfy $e_1>0$, $-c_1<0$ and $-c_2<0$.

For small enough $\tilde\delta$, in the scaled neighbourhoods
$B_{\tilde\delta}(\xi_2)$ the restriction of the system to the unstable
manifold of $\xi_2$ in the leading order is $\dot s=\mu s+\beta\bar s^2$,
where we have denoted $s=w_2+\ri q_2$.
 As in the proof of Theorem \ref{th1}, $(\rho_1,\theta_1)$ and
$(\rho_2,\theta_2)$ denote polar coordinates in
$H^{(out)}_1$ and $H^{(in)}_2$, respectively, such that
$v_1=\rho_1\cos\theta_1$, $q_1=\rho_1\sin\theta_1$,
$w_2=\rho_2\cos\theta_2$ and $q_2=\rho_2\sin\theta_2$.
We assume that the local bases near $\xi_1$ and $\xi_2$ are chosen in such
a way that the heteroclinic connection $\xi_2\to\xi_1$
goes along the directions $\arg(\theta_j)=0$ for both $j=1,2$, which implies $\beta>0$.
In the complement subspace the system is
approximated by the contractions $\dot u=-r_2 u$ and $\dot v=-c_2 v$.
In terms of the functions $\tau(r,\theta)$ and $\vartheta(r,\theta)$ introduced
in Lemma \ref{d3_30}, the map $\phi_2$ is
$$(v_2,q_2)=\phi_2(\rho_2,\theta_2)=
(v_{0,2}\re^{-c_2\tau(\rho_2,\theta_2)},\sin\vartheta(\rho_2,\theta_2)).$$

In physical space, there are two heteroclinic trajectories from $\xi_1$ to
$\xi_2$, with positive or negative $w_1$, and three trajectories from $\xi_2$
to $\xi_1$, with $\theta_2=0$, $2\pi/3$ or $4\pi/3$.
Let $H^{(out)}_1$ be the crossection of the heteroclinic trajectory $\xi_1\to\xi_2$
with positive $w_1$. We take a modified, compared to (\ref{glmap1}), map $\psi_1$:
\begin{equation}\label{glmap1n}
(\rho_2,\theta_2)=\psi_1(\rho_1,\theta_1)=\gamma(A\rho_1,\theta_1+\Theta)
\hbox{, where }\gamma=\kappa^l\rho^s.
\end{equation}
By choosing $s$, we obtain
\begin{equation}\label{glmap1nn}
(\rho_2,\theta_2)=\psi_1(\rho_1,\theta_1)=(A\rho_1,(-1)^l\Theta')
\hbox{, where }\Theta'=\Theta+2\pi s/3\hbox{ satisfies }-\pi/3<\Theta'<\pi/3.
\end{equation}

According to Lemma \ref{d3_30}(iv), for small $\rho_2$ and $\mu$
$$\re^{-c_2\tau(\rho_2,\theta_2)}\ll\sin\vartheta(\rho_2,\theta_2),$$
therefore we have $\psi_2(v_2,q_2)\approx(B_{12}q_2,B_{22}q_2)$.
For small $\theta_1$ we have $\sin\theta_1\approx\tan\theta_1\approx\theta_1$.
Taking into account (\ref{phmap1}), (\ref{glmap1nn}),
Lemma \ref{d3_1} and Lemma \ref{d3_30}(iii), we obtain that
\begin{equation}\label{mapg}
g(\rho_1,\theta_1)\approx
\biggl(C_1(\rho_1A\beta\sin^{1/3}3\tilde\Theta+\mu S(\tilde\Theta))^{3c_1/e_1},
C_2(\rho_1A\beta\sin^{1/3}3\tilde\Theta+\mu S(\tilde\Theta))\biggr),
\end{equation}
where $C_1=v_{0,1}3^{-1}\beta^{-3}|B_{12}|^{c_1/e_1}$ and
$C_2=(-1)^lv_{0,1}^{-1}3^{-1}\beta^{-3}B_{22}$.
(Here, we have denoted $\tilde\Theta=|\Theta'|$ and the
power $l$ in (\ref{glmap1nn}) is chosen in such a way that in
$(w_1,q_1)=\psi_2\phi_2\psi_1(\rho_1,0)$ the value of $w_1$ is positive
for small $\rho_1$.)

\medskip
(a) From (\ref{mapg}), the $\rho$-component of $g$ satisfies
$$g_{\rho}(\rho_1,\theta_1)> C_3\rho_1^{3c_1/e_1},
\hbox{ where }C_3=C_1(A\beta\sin^{1/3}\tilde\Theta)^{3c_1/e_1},$$
hence if $3c_1<e_1$ then for any $0<\delta<C_3^{e_1/(e_1-3c_1)}$ the iterates
$g^n(\rho_1,\theta_1)$ with initial $0<\rho_1<\delta$ satisfy
$g_{\rho}^n(\rho_1,\theta_1)>\delta$ for sufficiently large $n$.

\medskip
(b) Assume that $3c_1>e_1$.
Since $\rho_1=0$ is a solution to
\begin{equation}\label{mapr1}
g_{\rho}(\rho_1,\theta_1)-\rho_1=0
\end{equation}
for $\mu=0$, by the implicit function theorem, the equation (\ref{mapr1}) has
a unique solution $\rho_1(\mu)$, $0\le\mu\le\mu_1$, for some $\mu_1>0$ and
$\rd r_1/\rd\mu=3c_1C_1\mu^{(3c_1-e_1)/e_1}S^{3c_1/e_1}(\tilde\Theta)/e_1$ at $\mu=0$.
(In order to apply the implicit function theorem, we define
the function $h(\rho_1,\mu)=\rho_1-g_{\rho}(\rho_1,\theta_1)$ for
negative $\rho_1$ by setting $h(-\rho_1,\mu)=-h(\rho_1,\mu)$. Note, that
$g_{\rho}(\rho_1,\theta_1)$ does not depend on $\theta_1$.)
Therefore, for small $\mu$ the fixed point of the map $g$ (\ref{mapg}) can
be approximated by
$(\rho_p,\theta_p)=(C_1(\mu S(\tilde\Theta))^{3c_1/e_1},C_2\mu S(\tilde\Theta))$.
This fixed point is an intersection
of a periodic orbit with $H^{(out)}_1$. The distance from $(\rho_p,\theta_p)$
to $X$ depends on $\mu$ as $\mu^{3c_1/e_1}$, therefore the trajectory approaches
$X$ as $\mu\to0$. For a given small $\delta>0$, to find $\mu(\delta)$ we
approximate trajectories near $\kappa_1$ and $\kappa_2$ by linear (in
$(v_1, q_1)$ and $(v_2, q_2)$, respectively) maps, near $\xi_1$ we use
the approximation (\ref{lmap}) and near $\xi_2$ we employ Lemma \ref{d3_1}.
We do not go into details.

To study stability of the fixed point $(\rho_p,\theta_p)$, we consider the
difference $g(\rho_1+r',\theta_1+\theta')-g(\rho_1,\theta_1)$,
assuming that $r'$ and $\theta'$ are small and $(\rho_1,\theta_1)$ are
close to $(\rho_p,\theta_p)$. For the maps $\psi_1$, $\psi_2$ and $\phi_1$
we have
\begin{equation}\label{mappert}
\begin{array}{l}
\psi_1(\rho_1+r',\theta_1+\theta')-\psi_1(\rho_1,\theta_1)=(Ar',\theta'),\\
\psi_2(v_2+v',q_2+q')-\psi_2(v_2,q_2)=(B_{11}v'+B_{12}q',B_{21}v'+B_{22}q'),\\
\phi_1(w_1+w',q_1+q')-\phi_1(w_1,q_1)=
((v_{0,1}c_1/e_1)w_1^{c_1/e_1-1}w',(1+(q_1/v_{0,1})^2)^{-1}v_{0,1}^{-1}q').
\end{array}\end{equation}

Recall that (see Lemma \ref{d3_30}(v))
\begin{equation}\label{phi2r}
|\re^{-c_2(\tau(\rho_2+r',\theta_2+\theta'))}-
\re^{-c_2\tau(\rho_2,\theta_2)}|<\re^{-c_2\tau(\rho_2,\theta_2)}|
{r'\over\mu \rho_2}+3\theta' S(\theta_2)\tau(\rho_2,\theta_2)|.
\end{equation}
From (\ref{estt}), for small $\mu$ and $\rho_2\ll\mu$
\begin{equation}\label{phi2t}
|\vartheta(\rho_2+r',\theta_2+\theta')-\vartheta(\rho_2,\theta_2)|\approx
\vartheta^{2/3}(\rho_2,\theta_2)
(r'\sin^{1/3}3\tilde\Theta+\theta'\mu\sin^{-2/3}3\tilde\Theta/\beta).
\end{equation}
By arguments similar to the ones applied in the proof of Lemma \ref{d3_30}(iv),
for small $\mu$ and $\rho_2$ the r.h.s. of (\ref{phi2r}) is
asymptotically smaller than the r.h.s. of (\ref{phi2t}). Therefore, combining
(\ref{phi2t}) with (\ref{mappert}), we obtain that for small $\mu$
and $\rho_1\ll\mu$
$$g(\rho_1+r',\theta_1+\theta')-g(\rho_1,\theta_1)\approx
(A_1\mu^{(3c_1-e_1)/e_1}r'+A_2\mu^{3c_1/e_1}\theta',
A_3\mu^2r'+A_4\mu^3\theta'),$$
where the constants $A_j$ depend on $v_{0,1}$, $e_1$, $c_1$, $\beta$, $B_{12}$,
$B_{22}$ and $\tilde\Theta$.
Since $3e_1>c_1$, for small $\mu$ we have
$$|g(\rho_1+r',\theta_1+\theta')-g(\rho_1,\theta_1)|<
\mu^q|(r',\theta')|,$$
and $q=\min((3c_1-e_1)/(2e_1),1)$,
which implies that the bifurcating periodic orbit is asymptotically stable.
\qed

\begin{remark}\label{varorb}
In the proof of the Theorem we have considered the map
$g:\ H^{(out)}_1\to H^{(out)}_1$, $g=\phi_1\psi_2\phi_2\psi_1$, where
$\psi_1$ involves the symmetry $\gamma=\kappa^l\rho^s$ and the choice of $l$ and $s$ depends on
$\Theta$ and $B_{12}$. Hence, depending on the values of these constants,
there can exist geometrically different periodic orbits in the vicinity
of the group orbit of the heteroclinic cycle $\xi_1\to\xi_2\to\xi_1$
with different number of symmetric copies of heteroclinic connections
$\kappa_1$ and $\kappa_2$. The cycle shown in Fig.~\ref{fig2} involves two
symmetric copies of the connection $\kappa_1$ and one connection $\kappa_2$.
For a different group $\Gamma\subset$\,SO(4) considered in Section
\ref{sec:otherexample}, we present examples of various periodic orbits which
are obtained by varying the coefficients of the respective normal form, see Fig.~\ref{fig4}.
\end{remark}

\subsection{A numerical example} \label{sec:num_Gamma}

In this subsection we present an example of a $\Gamma$-equivariant vector
field possessing a heteroclinic cycle, introduced in proposition
\ref{prop:admissible}, where the expanding eigenvalue of $df(\xi_2)$ is
small. In agreement with Theorem \ref{thperorb}, asymptotically stable
periodic orbits exist in the vicinity of the cycle.
To construct the numerical example, we start from a Lemma that determines the
structure of $\Gamma$-equivariant vector fields (the proof is left to the reader):
\begin{lemma} \label{lem:cubicGammaequiv}
Any $\Gamma$-equivariant vector field is
\begin{equation}\label{eq:equivariantany}
\begin{array}{l}
\dot z_1= A_{n_1n_2n_3n_4}z_1^{n_1}\bar z_1^{n_2}z_2^{n_3}\bar z_2^{n_4}\\
\dot z_2= B_{n_1n_2n_3n_4}z_1^{n_1}\bar z_1^{n_2}z_2^{n_3}\bar z_2^{n_4},\\
\end{array}
\end{equation}
where $A_{n_1n_2n_3n_4}$ and $B_{n_1n_2n_3n_4}$ are real, $A_{n_1n_2n_3n_4}\ne0$
whenever $n_3-n_4=0(\mmod\,3)$ and $B_{n_1n_2n_3n_4}\ne0$ whenever $n_3-n_4=1(\mmod\,3)$.
\end{lemma}

Keeping in (\ref{eq:equivariantany}) all terms of order one and two, and
several terms of the third order, we consider the following vector field:
\begin{equation}\label{eq:equivariantcubic}
\begin{array}{l}
\dot z_1= a_1z_1+a_2\bar z_1+a_3z_1^2+a_4\bar z_1^2+a_5z_1\bar z_1+
a_6z_2\bar z_2+a_7z_1^2\bar z_1+a_8z_1z_2\bar z_2+a_9z_2^3 +a_{10}\bar z_2^3 \\
\dot z_2=b_1z_2+b_2z_1z_2+b_3\bar z_1z_2+b_4\bar z_2^2+b_5z_2z_1\bar z_1+b_6z_2^2\bar z_2,
\end{array}
\end{equation}
where
\begin{equation}\label{assume}
\hbox{we set }a_3+a_4+a_5=0,\ a_7=a_8=-1,
\hbox{ and denote }\alpha=a_1+a_2,\ \alpha'=a_1-a_2.
\end{equation}
The system (\ref{eq:equivariantany}) restricted into $L=\Fix\D_3$ is
$$
\dot x_1=\alpha x_1-x_1^3,
$$
which implies that the steady states $\xi_1$ with $x_1=\alpha^{1/2}$
and $\xi_2=-\xi_1$ exist whenever $\alpha>0$. The radial eigenvalues
of $\xi_j$ are $-2\alpha$. The
non-radial eigenvalues of $df(\xi_1)$ are
\begin{equation}\label{evs-xi1}
\begin{array}{l}
-c_1= a_1-a_2+ 2(a_3-a_4)\sqrt{\alpha} - \alpha \mbox{~along the eigendirection~} y_1 \\
e_1= b_1 + (b_2+b_3)\sqrt{\alpha} + b_5\alpha \mbox{~along the eigendirections~} x_2\hbox{ and }y_2
\end{array}
\end{equation}
and the eigenvalues at $df(\xi_2)$ are
\begin{equation}\label{eigenvalues}
\begin{array}{l}
e_2= a_1-a_2-2(a_3-a_4)\sqrt{\alpha} - \alpha \mbox{~along the eigendirection~} y_1 \\
-c_2= b_1 - (b_2+b_3)\sqrt{\alpha} + b_5\alpha \mbox{~along the eigendirections~} x_2\hbox{ and }y_2
\end{array}
\end{equation}
The eigenvalues $e_1$ and $-c_2$ are double.

With (\ref{assume}) taken into account, the restriction of the system
(\ref{eq:equivariantcubic}) into the plane $P_1$ is
\begin{equation}\label{eq:P1}
\begin{array}{l}
\dot x_1= \alpha x_1 - 2(a_3+a_4)y_1^2 - x_1(x_1^2+y_1^2) -x_1x_2^2 \\
\dot y_1= \alpha' y_1 + 2(a_3-a_4)x_1y_1 - y_1(x_1^2+y_1^2) \\
\end{array}
\end{equation}
and the one into the plane $P_2$ is
\begin{equation}\label{eq:P2}
\begin{array}{l}
\dot x_1= \alpha x_1 + a_6x_2^2 - x_1^3 -x_1x_2^2 + (a_9+a_{10})x_2^3 \\
\dot x_2 = b_1 x_2 + (b_2+b_3)x_1x_2 + b_4 x_2^2 + b_5 x_1^2x_2 + b_6 x_2^3
\end{array}
\end{equation}
In both planes the system has a form which is well-known to produce saddle-sink
connections between equilibria $\xi_1$ and $\xi_2$
if certain conditions are fulfilled, see e.g. \cite{cho01}. A sufficient
condition for existence of a heteroclinic connection $\xi_1\to\xi_2$ in $P_1$ is
\begin{equation}\label{conP1}
a_3+a_4>0,\ a_3-a_4>0,\
a_3-a_4 +((a_3-a_4)^2+\alpha')^{1/2}<\alpha^{1/2}.
\end{equation}
The expression for a condition of the existence of a connection $\xi_2\to\xi_1$
in $P_2$ is too bulky, and therefore is not presented.

\begin{figure}[p]

\vspace*{-7cm}
\hspace*{-2cm}\includegraphics[width=16cm]{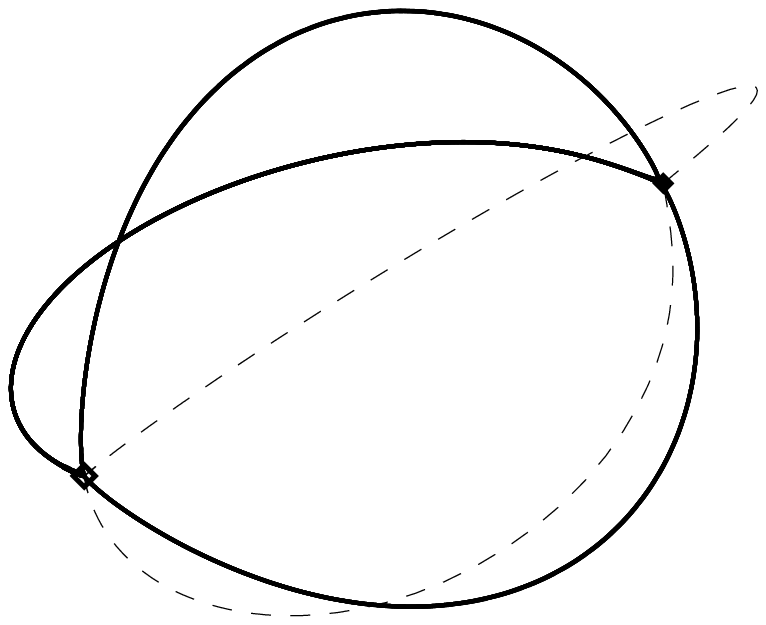}

\vspace*{-62mm}
\hspace*{113mm}{\large (a)}

\vspace*{9mm}
\hspace*{-2cm}\includegraphics[width=19cm]{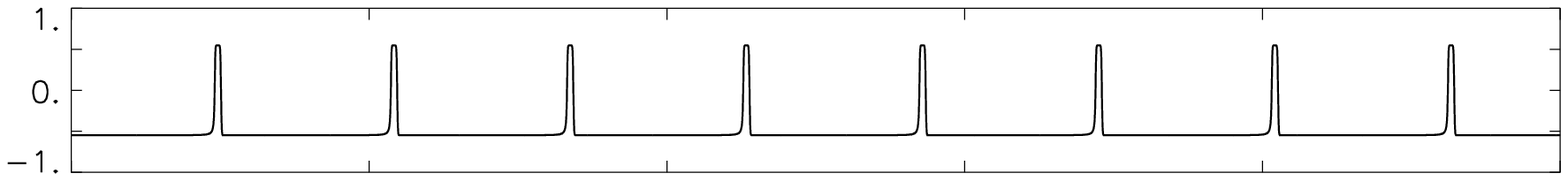}

\vspace*{-18mm}
\hspace*{-2cm}\includegraphics[width=19cm]{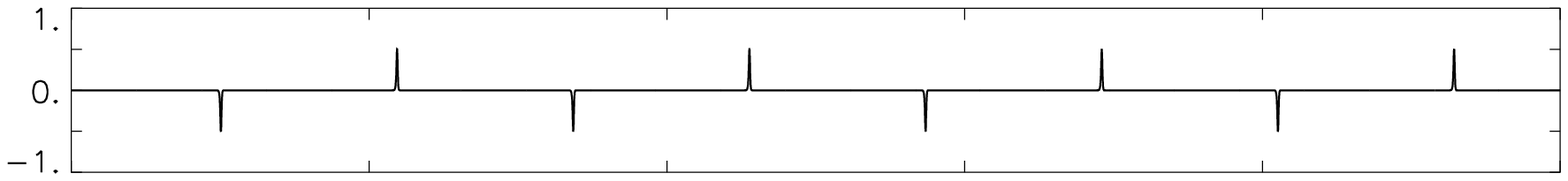}

\vspace*{-12mm}
\hspace*{-3mm}{10000\hspace*{22mm}12000\hspace*{22mm}14000\hspace*{22mm}16000\hspace*{22mm}18000\hspace*{22mm}20000}

\vspace*{1mm}
\hspace*{152mm}{\large $t$}

\vspace*{1mm}
\hspace*{113mm}{\large (b)}

\vspace*{-53mm}
\hspace*{-12mm}{\large $x_1$}

\vspace*{15mm}
\hspace*{-12mm}{\large $y_1$}

\vspace*{34mm}
\noindent
\caption{Projection of the periodic orbits (solid line) and of the
heteroclinic connections $\xi_2\to \xi_1$ (dashed lines) into the plane
$<{\bf v}_1,{\bf v}_2>$, where ${\bf v}_1=(4,2,4,1.5)$ and
${\bf v}_2=(2,4,-1.5,4)$, (a). The steady state $\xi_1$ is denoted by
a hollow circle and $\xi_2$ by filled one.
Temporal behaviour of $x_1$ and $y_1$, (b).}
\label{fig2}\end{figure}

We choose the values of the coefficients so that the conditions
for existence of heteroclinic connections in $P_1$ and $P_2$ are satisfied,
with $\xi_1$ being a saddle in $P_1$ and a sink in $P_2$ and vice versa for
$\xi_2$. We also choose the coefficients such that the stability condition
$3c_1>e_1$ (see Theorem \ref{thperorb}) is satisfied.\\
Finally we adjust coefficients so that the unstable, double eigenvalue $e_2$
at $\xi_2$ verifies is small.
The simulations are performed with the following values of coefficients:
\begin{equation}
\begin{array}{l}\label{coefficients}
\alpha=0.3,~\alpha'=0.2,~a_3=0.3,~a_4=-0.05,~a_5=-0.25,~a_6=0.6, \\
a_7=a_8=-1,~a_9=0.1,~a_{10}=0.15; \\
b_1=0.2,~b_2=-0.1,~b_3=-0.09,~b_4=-0.1,~b_5=-1,~b_6=-1,
\end{array}
\end{equation}
which implies that the eigenvalues are
$-c_1=-0.2041$, $e_1=0.2834$, $-c_2=-0.4834$, $e_2=0.0041$. \\
In Fig.~\ref{fig2} we show projection of the periodic orbit and of
the group orbit of heteroclinic connections, comprising the cycle, on a plane
in $\R^4$, and time series of $x_1$ and $y_1$. The
periodic orbit follows two heteroclinic connections $\xi_1\to\xi_2$ in $P_1$,
and only one $\xi_2\to\xi_1$ in $P_2$, see remark \ref{varorb}.

In Section \ref{sec:otherexample} we shall explore numerically the dynamics
near the pseudo-simple heteroclinic cycle with a more involved symmetry subgroup of SO(4).

\section{Stability of pseudo-simple heteroclinic cycles when
$\Gamma\subset$\,O(4), but $\Gamma\not\subset$\,SO(4)}\label{sec:nonSO(4)}

Here we show that presence of a reflection in the group $\Gamma$ can
completely change the asymptotic dynamics of the system (\ref{eq_ode}) near
a pseudo-simple heteroclinic cycle.
As in the previous Section we hold on an example, which is built as
an extension of the group $\Gamma$ generated by transformations \eqref{gen:Gamma}.
Let's introduce the reflection
\begin{equation}
\sigma:~(z_1,z_2)\mapsto (z_1,\bar z_2)
\end{equation}
and define $\tilde\Gamma=\Gamma\cup\sigma\Gamma$. This group admits the same
one-dimensional and two-dimensional fixed point subspaces as $\Gamma$, therefore it
also admits (similar) pseudo-simple heteroclinic cycles.
However, it also has two types of three-dimensional fixed point subspaces: \\
$V=\{(x_1,y_1,x_2,0)\}=\Fix(\sigma)$, \\
$W=\{(x_1,0,x_2,y_2)\}=\Fix(\sigma\kappa)$. \\
Note that (i) $V=P_1\oplus P_2$ (the invariant planes introduces in Section
\ref{exa1}) and (ii) $V$ has two symmetric copies $V'=\rho V$ and
$V''=\rho^2 V$, while $W$ is a singleton. It follows from point (i) that if it exists, the heteroclinic cycle lies entirely in $V$.

\subsection{A Theorem}

\begin{theorem}\label{thsect5}
Consider the $\tilde\Gamma$-equivariant system
\begin{equation}\label{eqal_ode5}
\dot{\bf x}=f({\bf x}),\hbox{ where }
f(\gamma{\bf x})=\gamma f({\bf x})\mbox{ for all }\gamma\in\tilde\Gamma
\end{equation}
and $f$ is a smooth map. Suppose that the system possesses equilibria $\xi_1$ and $\xi_2$ and
a pseudo-simple heteroclinic cycle, introduced in proposition
\ref{prop:admissible}, where the heteroclinic connection $\xi_2\to\xi_1$ belongs
to the half-plane of $P_2$ satisfying $x_2>0$. We assume that
$\partial^2 f_{x_2}/\partial x_2^2>0$, where $f_{x_2}$ denotes the
$x_2$-component of $f$. Let $-c_j$ and $e_j$ be the stable and
unstable eigenvalues of $df(\xi_j)$, respectively,
such that $-c_1$ and $e_2$ have multiplicity 2.\\
Then, for sufficiently small expanding eigenvalue $e_2$ the cycle
is f.a.s. (see definition \ref{def:fragmstable}).
\end{theorem}

We begin with a proof of a Lemma about properties of trajectories
in a generic $\D_3$-equivariant second order dynamical system (\ref{sysA}),
where we assume that $\alpha>0$ and $\beta>0$. Recall, that in polar
coordinates the system takes the form (\ref{pcor1}),
the system is considered in the sector $0\le\theta\le\pi/3$, by
$\tau(r_0,\theta_0)$ we have denoted the time it takes the trajectory of the
system (\ref{pcor1}) starting at $(r_0,\theta_0)$ to reach $r=1$ and
by $\vartheta(r_0,\theta_0)$ the value of $\theta$ at $r=1$.

\begin{lemma}\label{d5}
If $0<\theta_0<\pi/6$ then $\vartheta(r_0,\theta_0)$ satisfies
\begin{equation}\label{est51}
\sin3\vartheta(r_0,\theta_0)<\biggl({r_0\beta\cos3\theta_0+\alpha\over
\beta\cos3\theta_0+\alpha}\biggr)^3\sin3\theta_0.
\end{equation}
\end{lemma}

\proof
From (\ref{pcor1}),
$$
\frac{d}{dt}(r^6\sin^23\theta) = 6\alpha r^6\sin^23\theta,
$$
which implies that
$$
r^6(t)\sin^23\theta(t) = e^{6\alpha t}r_0^6\sin^23\theta_0.
$$
Setting $t=\tau(r_0,\theta_0)$ we obtain that
\begin{equation}\label{sol:sin3theta(tau)}
\sin3\theta(\tau) = e^{3\alpha\tau}r_0^3\sin3\theta_0.
\end{equation}
As well, (\ref{pcor1}) implies that $\dot\theta<0$ and that for
$\theta_1<\theta_2$ we have $\dot r(r,\theta_1)>\dot r(r,\theta_2)$.
Hence, the trajectory $(r(t),\theta(t))$
starting at $r(0)=r_0$ and $\theta(0)=\theta_0$ satisfies
$$\dot r(r(t),\theta(t))>\dot r(r(t),\theta_0)=\alpha r+\beta\cos3\theta_0r^2
\hbox{ for all }t>0.$$
Therefore, by the same arguments as we employed in the proof of Lemma \ref{d3_30}(i,ii)
\begin{equation}\label{estl5}
{\rm e}^{\alpha\tau(r_0,\theta_0)}<{r_0\beta\cos3\theta_0+\alpha\over
r_0(\beta\cos3\theta_0+\alpha)}.
\end{equation}
Combining (\ref{sol:sin3theta(tau)}) and (\ref{estl5}) we obtain statement
of the Lemma.
\qed

\noindent
{\bf Proof of Theorem \ref{thsect5}:}
We aim at finding conditions such that the first return map
$g=\psi_1\phi_1\psi_2\phi_2:\ H^{(in)}_2\to H^{(in)}_2$ is a contraction and
hence converges to $0$ as the number of iterates tends to infinity.
(Note, that here we consider the map $H^{(in)}_2\to H^{(in)}_2$, while
in Theorems \ref{th1} and \ref{thperorb} it was $H^{(out)}_1\to H^{(out)}_1$.)
The local and global maps are defined the same way, as it was in the proofs
of Theorems \ref{th1} and \ref{thperorb}. Because of the reflection $\sigma$,
the constants $\Theta$, $B_{12}$ and $B_{21}$ in (\ref{glmap1}) and (\ref{psi2})
vanish, hence global maps take the form
\begin{equation}\label{glmap12}
(\rho_2,\theta_2)=\psi_1(\rho_1,\theta_1)=(A\rho_1,\theta_1),\quad
(w_1,q_1)=\psi_2(v_2,q_2)=(B_{11}v_2,B_{22}q_2).
\end{equation}
Recall, that $(\rho_1,\theta_1)$ and $(\rho_2,\theta_2)$ denote polar
coordinates in $H^{(out)}_1$ and $H^{(in)}_2$, respectively, such that
$v_1=\rho_1\cos\theta_1$, $q_1=\rho_1\sin\theta_1$,
$w_2=\rho_2\cos\theta_2$ and $q_2=\rho_2\sin\theta_2$.
The map $\phi_1$ is given by (\ref{phmap1}).
The map $\phi_2$ is
$$(v_2,q_2)=\phi_2(\rho_2,\theta_2)=
(v_{0,2}\re^{-c_2\tau(\rho_2,\theta_2)},\tan\vartheta(\rho_2,\theta_2)).$$
The composition $g=\psi_1\phi_1\psi_2\phi_2$ therefore is
\begin{equation} \label{eq:g}
g(\rho_2,\theta_2) =\left(Av_{0,1}\biggl(B_{11}v_{0,2}
{\rm e}^{-c_2\tau(\rho_2,\theta_2)}\biggr)^{c_1/e_1},
\arcsin(B_{22}\tan\vartheta(\rho_2,\theta_2) \right).
\end{equation}

The condition $\partial^2 f_{x_2}/\partial x_2^2>0$ implies that the
restriction of the system (\ref{eqal_ode5}) into the unstable manifold
of $\xi_2$, and which is of the form (\ref{sysA}), satisfies $\beta>0$.
The condition $\alpha>0$ is evidently satisfied because of the instability
of $\xi_2$, hence we can apply Lemmas \ref{d3_30} and \ref{d5} to estimate
$\tau(\rho_2,\theta_2)$ and $\vartheta(\rho_2,\theta_2)$.
Denote by $g_{\rho}$ and $g_{\theta}$ the first and second components of $g$
in \eqref{eq:g}. From Lemmas \ref{d3_30}(i,ii) and \ref{d5} we can write the
estimates
\begin{equation}\label{eq:g_r,g_theta}
\begin{array}{l}
g_{\rho}(\rho_2,\theta_2) \leq C\left(
{\rho_2(\beta+e_2)\over \rho_2\beta+e_2}\right)^{c_1c_2/e_1e_2} \\
g_\theta(r_2,\theta_2)\leq\arcsin\left(B_{22}\tan\left(1/3\arcsin
\left[\left(\frac{\rho_2\beta\cos3\theta_2+e_2}{\beta\cos3\theta_2+e_2}\right)^3\cdot\sin(3\theta_0) \right]\right)\right)
\end{array}
\end{equation}
where $C=Av_{0,1}(B_{11}v_{0,2})^{c_1/e1}$.
We show that $g_{\rho}$ and $g_\theta$ are contractions under the conditions
stated in the Theorem. \\
The function $g_{\rho}(\rho_2,\theta_2)$ is a smooth function of $\rho_2$
for all $\rho_2>0$. Its derivative is positive and has the expression
$$
g_{\rho}'(\rho_2)=hC\frac{e_2(\beta+e_2)^h}
{(\rho_2\beta+e_2)^{h+1}}\rho_2^{h-1},
$$
where we have denoted $h=\frac{c_1c_2}{e_1e_2}$.
The condition $h>1$ implies that
$\lim\limits_{\substack{\rho \to 0}}g_{\rho}'(\rho)=0$.
Therefore, there exists $\eta>0$ such that $g_{\rho}'(\rho_2)<1$ and
$g_{\rho}$ is a contraction when $\rho_2<\eta$. \\
To prove that $g_\theta$ is a contraction let us assume that $\theta_2$ is
chosen small enough to allow for the approximation
$\tan 3\theta_2\approx\sin 3\theta_2\approx 3\theta_2$ (this is not necessary but simplifies the calculations). Then
$$
g_\theta(\theta_2)\leq B_{22}\left(\frac{\rho_2\beta\cos3\theta_2+e_2}
{\beta\cos3\theta_2+e_2}\right)^3\theta_2.
$$
If we chose $e_2$ and $\rho_2$ small enough ($\beta$ is a fixed positive
constant), then the factor in front of $\theta_2$ is smaller than 1.
Therefore, if $\rho_2$ converges to 0 by iterations of $g$, the same is true
for $\theta_2$. Combining the conditions for $g_{\rho}$ and $g_\theta$ we get the result.
\qed

Remark that the convergence to 0 of $\theta_2$ is slow compared to that of
$\rho_2$.

\begin{remark}
Suppose that, similarly to Theorem \ref{thperorb}, we consider a
$\tilde\Gamma$-equivariant system which depends on a parameter $\mu$, possesses
a heteroclinic cycle for $0<\mu<\mu_0$ and $e_2=0$ for $\mu=0$.
Then by Theorem \ref{thsect5} there exists $\mu_1>0$ such that the heteroclinic cycle is f.a.s. for
$0<\mu<\mu_1$. Moreover, the condition
$\partial^2 f_{x_2}/\partial x_2^2>0$ (which is equivalent to $\beta>0$)
is always satisfied in such a system, by the reasons given in the
proof of Theorem \ref{thperorb}. Note however, that closeness to the
bifurcation point, and even smallness of $e_2$, are not necessary
for fragmentary asymptotic stability of the cycle.
\end{remark}

\subsection{A numerical illustration of Theorem \ref{thsect5}} \label{subsec:vectorfieldO(4)}

The equivariant structure of $\tilde\Gamma$-equivariant vector fields is easily deduced from that of $\Gamma$ equivariant vector fields.
\begin{lemma} \label{lem:cubictildeGammaequiv}
Any $\tilde\Gamma$-equivariant vector field has the form
\eqref{eq:equivariantany} (see Lemma \ref{lem:cubicGammaequiv}) with the following
additional conditions on the coefficients: $A_{n_1n_2n_3n_4}=A_{n_1n_2n_4n_3}$
and $B_{n_1n_2n_3n_4}=B_{n_2n_1n_3n_4}$.
\end{lemma}
In this subsection we consider the same system (\ref{eq:equivariantany}) as
we considered in subsection \ref{sec:num_Gamma}. Because of the presence of
the symmetry $\sigma$ (see Lemma \ref{lem:cubictildeGammaequiv}), coefficients
of this system satisfy $a_9=a_{10}$ and $b_2=b_3$. As in subsection
\ref{sec:num_Gamma}, we denote $\alpha=a_1+a_2$, $\alpha'=a_1-a_2$
and assume $a_3+a_4+a_5=0$ and $a_7=a_8=-1$. Therefore, the system is
\begin{equation}\label{exasec5}
\begin{array}{l}
\dot z_1= a_1z_1+a_2\bar z_1+a_3z_1^2+a_4\bar z_1^2+a_5z_1\bar z_1+
a_6z_2\bar z_2-z_1^2\bar z_1-z_1z_2\bar z_2+a_9(z_2^3+\bar z_2^3) \\
\dot z_2=b_1z_2+b_2(z_1z_2+\bar z_1z_2)+b_4\bar z_2^2+b_5z_2z_1\bar z_1+b_6z_2^2\bar z_2,
\end{array}
\end{equation}

The conditions for the existence of the heteroclinic cycle in the invariant
space $V$ are similar to those given in Section \ref{sec:num_Gamma}.
It is however interesting to analyse the dynamics in the invariant space $W$.
Consider the restriction of system (\ref{exasec5}) into the subspace $W$:
\begin{equation}\label{eq:W}
\begin{array}{l}
\dot x_1= \alpha x_1+a_3(x_2^2+y_2^2)-x_1^3-x_1(x_2^2+y_2^2)+2a_9(x_2^3-3x_2y_2^2) \\
\dot x_2 =b_1x_2+2b_2x_1x_2+b_4(x_2^2-y_2^2)+b_5x_2x_1^2+b_6x_2(x_2^2+y_2^2) \\
\dot y_2= b_1y_2+2b_2x_1y_2-2b_4x_2y_2+b_5y_2x_1^2+b_6y_2(x_2^2+y_2^2).\\
\end{array}
\end{equation}
In the limit $b_4=0$ any plane containing the axis $L$ is flow-invariant by
\eqref{eq:W}. In fact, any such plane is a copy of $P_2$ by some rotation around
$L$. Therefore, whenever a saddle-sink connection exists in $P_2$, it generates
a two-dimensional manifold of saddle-sink connections in $W$ (an "ellipsoid" of connections).
Small $\tilde\Gamma$-equivariant perturbations of the system, in particular
switching on $b_4\neq 0$, do not destroy this manifold. If in addition
the conditions (\ref{conP1}) are satisfied, then
a continuum of robust heteroclinic cycles connecting $\xi_1$ and $\xi_2$ exists for the $\tilde\Gamma$-equivariant vector field.

We take the same values of the coefficients as in \eqref{coefficients}, except
that $a_9=a_{10}=0.1$ and $b_2=b_3=-0.6$ so that the system is
$\tilde\Gamma$-equivariant.

\begin{figure}[p]

\vspace*{-7cm}
\hspace*{-2cm}\includegraphics[width=16cm]{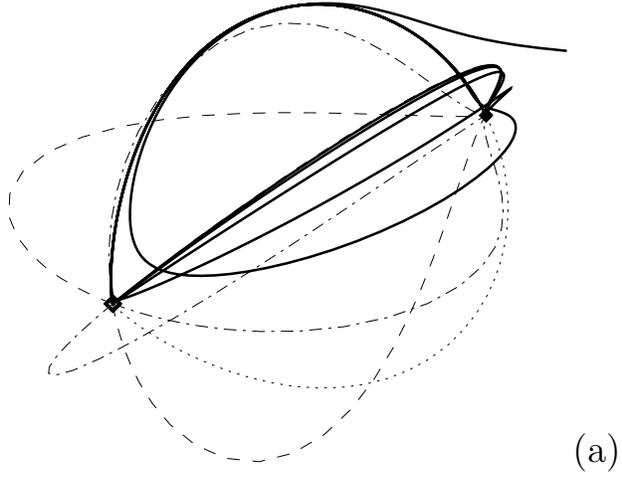}

\vspace*{-52mm}
\hspace*{103mm}{\large (a)}

\vspace*{9mm}
\hspace*{-2cm}\includegraphics[width=9.5cm]{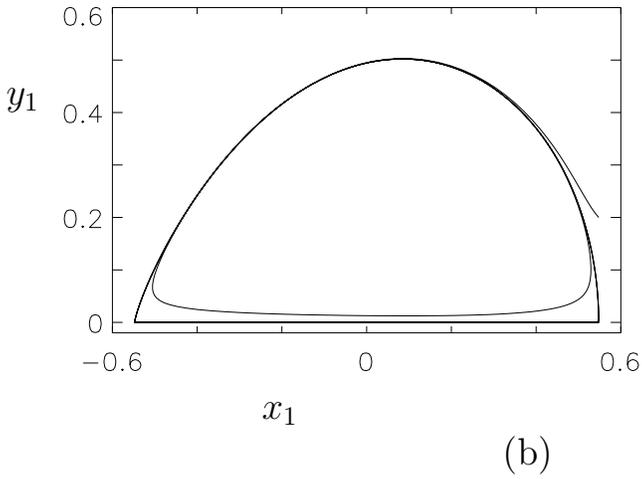}
\hspace*{-5mm}\includegraphics[width=9.5cm]{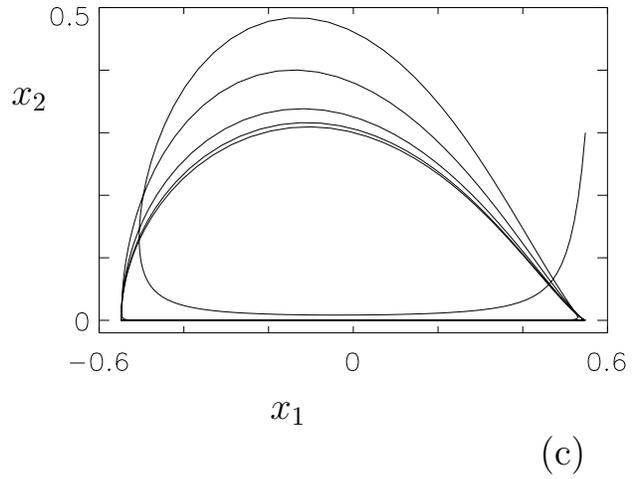}

\vspace*{-7mm}
\hspace*{21mm}{\large $x_1$\hspace*{89mm}$x_1$}

\vspace*{1mm}
\hspace*{53mm}{\large (b)\hspace*{91mm}(c)}

\vspace*{-53mm}
\hspace*{-13mm}{\large $y_1$\hspace*{89mm}$x_2$}

\vspace*{54mm}
\noindent
\caption{Projection of a trajectory of the system (\ref{exasec5}) converging
to a heteroclinic cycle (solid line) and of
heteroclinic connections (dotted, dashed and dashed-dotted lines, different
types of lines label different group orbits of heteroclinic connections) into
the plane $<{\bf v}_1,{\bf v}_2>$, where ${\bf v}_1=(4,-2,4,1.5)$ and
${\bf v}_2=(2,4,-1.5,4)$, (a). Projection of this trajectory into the planes
$(x_1,y_1)$ (b) and $(x_1,x_2)$ (c).}
\label{fig3}\end{figure}

Note, that when the symmetry $\sigma$ is absent, the heteroclinic cycle is
completely unstable but a nearby asymptotically stable periodic orbit is observed
exists, in accordance with the results of previous Sections. When the symmetry
$\sigma$ is present, the trajectory starting from some point in a
neighborhood of $\xi_2$ asymptotically converges to the heteroclinic cycle
(or to its symmetric copy) as $t\rightarrow\infty$, see Fig.~\ref{fig3}.
The convergence in the $(x_1,x_2)$ plane is slower than in the $(x_1,y_1)$
plane, because the convergence $\theta\to0$ is slower than $\rho\to0$.

\section{Another example with a group $\Gamma\subset SO(4)$}\label{sec:otherexample}

The groups which have been considered in Sections \ref{sec4} and
\ref{sec:nonSO(4)} had the advantage of allowing for an elementary algebraic
description, which gives an opportunity to study in detail the structure
on invariant subspaces and related dynamics.
In this Section we aim at providing an example of an admissible group with
pseudo-simple heteroclinic cycles, with a non-elementary group structure.
In principle, all such groups can be listed as we did for admissible groups
with simple heteroclinic cycles in \cite{pc13}. \\
We then build an equivariant system for this group and study numerically its
asymptotic behavior near the pseudo-simple heteroclinic cycle.
We begin with a description of the bi-quaternionic presentation of finite
subgroups of $SO(4)$ together with some geometrical properties of these actions,
which we employ to identify invariant subspaces.

\subsection{Subgroups of $SO(4)$: presentation and notation}\label{sec:quaternions}
Here we briefly introduce the bi-quaternionic presentation and notation for finite subgroups of $SO(4)$, which has been used in \cite{pc13} to determine all admissible groups for simple heteroclinic cycles in $\R^4$. See \cite{pdv} for more details.

A real quaternion is a set of four real numbers, ${\bf q}=(q_1,q_2,q_3,q_4)$.
Multiplication of quaternions is defined by the relation
\begin{equation}\label{mqua}
\begin{array}{ccc}
{\bf q}{\bf w}&=&(q_1w_1-q_2w_2-q_3w_3-q_4w_4,q_1w_2+q_2w_1+q_3w_4-q_4w_3,\\
&&q_1w_3-q_2w_4+q_3w_1-q_4w_2,q_1w_4+q_2w_3-q_3w_2+q_4w_1).
\end{array}\end{equation}
$\tilde{\bf q}=(q_1,-q_2,-q_3,-q_4)$ is the conjugate of $\bf q$, and
${\bf q}\tilde{\bf q}=q_1^2+q_2^2+q_3^2+q_4^2=|{\bf q}|^2$ is the square
of the norm of $\bf q$. Hence $\tilde{\bf q}$ is also the inverse ${\bf q}^{-1}$ of a unit
quaternion $\bf q$. We denote by $\mH$ the multiplicative group of unit
quaternions; obviously, its unity element is $(1,0,0,0)$.

Due to existence of a 2-to-1 homomorphism of $\mH$ onto SO(3) (see \cite{pdv}),
finite subgroups of $\mH$ are named after respective subgroups of SO(3). They are:
\begin{equation}\label{finsg}
\renewcommand{\arraystretch}{1.5}
\begin{array}{ccl}
\Z_n&=&\displaystyle{\oplus_{r=0}^{n-1}}(\cos2r\pi/n,0,0,\sin2r\pi/n)~,~~\D_n=\Z_{2n}\oplus\displaystyle{\oplus_{r=0}^{2n-1}}(0,\cos r\pi/n,\sin r\pi/n,0)\\
\V&=&((\pm1,0,0,0))~,~~\T=\V\oplus(\pm{1\over2},\pm{1\over2},\pm{1\over2},\pm{1\over2})\\
\mO&=&\T\oplus\sqrt{1\over2}((\pm1,\pm1,0,0))~,~~\I=\T\oplus{1\over2}((\pm\tau,\pm1,\pm\tau^{-1},0))\\
\end{array}\end{equation}
where $\tau=(\sqrt{5}+1)/2$. Double parenthesis denote all possible permutations
of quantities within the parenthesis and for $\I$ only even permutations of
$(\pm\tau,\pm1,\pm\tau^{-1},0)$ are elements of the group.

\medskip
By regarding the four numbers $(q_1,q_2,q_3,q_4)$ as Euclidean coordinates
of a point in $\R^4$, any pair of unit quaternions $({\bf l},{\bf r})$ can
be related to the transformation ${\bf q}\to{\bf lqr}^{-1}$, which is
an element of the group SO(4). The respective
mapping $\Phi:\mH\times\mH\to$\,SO(4) is a homomorphism onto,
whose kernel consists of two elements, $(1,1)$ and $(-1,-1)$.
Therefore, a finite subgroup of SO(4) is a subgroup of a product of two
finite subgroups of $\mH$. There is however an additional subtlety.
Let $\Gamma$ be a finite subgroup of SO(4) and ${\cal G}=\Phi^{-1}(\Gamma)=\{({\bf l}_j,{\bf r}_j),~1\le j\le J\}$.
Denote by $\bf L$ and $\bf R$ the finite subgroups of $\mH$
generated by ${\bf l}_j$ and ${\bf r}_j$, $1\le j\le J$, respectively.
To any element ${\bf l}\in \bf L$ there are several corresponding elements
${\bf r}_i$, such that $({\bf l},{\bf r}_i)\in\mH$, and similarly for any
${\bf r}\in{\bf R}$. This establishes a correspondence between $\bf L$ and
$\bf R$. The subgroups of $\bf L$ and $\bf R$ corresponding to
the unit elements in $\bf R$ and $\bf L$ are denoted by
${\bf L}_K$ and ${\bf R}_K$, respectively, and the group ${\cal G}$ by
$({\bf L}\rl{\bf L}_K;{\bf R}\rl{\bf R}_K)$.

\begin{lemma} (see proof in \cite{op13})\label{lem2}
Let $N_1$ and $N_2$ be two planes in $\R^4$ and $p_j$, $j=1,2$, be
the elements of SO(4) which act on $N_j$ as identity, and on $N_j^{\perp}$ as $-I$,
and $\Phi^{-1}p_j=({\bf l}_j,{\bf r}_j)$, where $\Phi$ is the homomorphism
defined above. Denote by $({\bf l}_1{\bf l}_2)_1$
and $({\bf r}_1{\bf r}_2)_1$ the first components of the respective quaternion
products. The planes $N_1$ and $N_2$ intersect if and only if
$({\bf l}_1{\bf l}_2)_1=({\bf r}_1{\bf r}_2)_1=\cos\alpha$ and $\alpha$ is
the angle between the planes.
\end{lemma}

\begin{lemma}\label{lem3}
Consider $g\in$\,SO(4), $\Phi^{-1}g=((\cos\alpha,\sin\alpha{\bf v});(\cos\beta,\sin\beta{\bf w}))$.\\
Then $\dim\Fix<g>=2$ if and only if $\cos\alpha=\cos\beta$.
\end{lemma}

\begin{lemma}\label{lem4}
Consider $g,s\in$\,SO(4), where
$\Phi^{-1}g=((\cos\alpha,\sin\alpha{\bf v});(\cos\alpha,\sin\alpha{\bf w}))$ and\\
$\Phi^{-1}s=((0,{\bf v});(0,{\bf w}))$.\\
Then $\Fix<g>=\Fix<s>$.
\end{lemma}

\subsection{The group $\Gamma=(\D_3\rl\Z_2;\mO\rl\V)$}\label{d3v}

The group $(\D_3\rl\Z_2;\mO\rl\V)$ was proposed in \cite{pc13} as an example of
a subgroup of O(4) admitting pseudo-simple heteroclinic cycles.
It is the (unique) four dimensional irreducible representation of the group
$GL(2,3)$ ($2\times 2$ invertible matrices over the field $\Z_3$). This group
is generated by the elements $\sigma_1$ (order 8) and $\sigma_2$ (order 2) below:
\begin{equation*}
\sigma_1 = \left(\begin{array}{cc}0&2\\2&2\end{array}\right),
~\sigma_2 = \left(\begin{array}{cc}2&0\\0&1\end{array}\right)
\end{equation*}

In quaternionic form its elements are:
\begin{eqnarray*}\label{listex1}
((1,0,0,0)&;&\V),\\
((\cos(\pi/3),0,0,\sin(\pi/3))&;&(1,1,1,1)/2~~\V),\\
((\cos(2\pi/3),0,0,\sin(2\pi/3))&;&(-1,1,1,1)/2~~\V),\\
((0,1,0,0)&;&(1,1,0,0)/\sqrt{2}~~\V),\\
((0,\cos(\pi/3),\sin(\pi/3),0)&;&(1,0,0,1)/\sqrt{2}~~\V),\\
((0,\cos(2\pi/3),\sin(2\pi/3),0)&;&(1,0,1,0)/\sqrt{2}~~\V).
\end{eqnarray*}

Lemmas \ref{lem3} and \ref{lem4} imply that the group has two isotropy types of
subgroups $\Sigma$ satisfying $\dim\Fix\Sigma=2$.
They are generated either by an order three element,
$$\epsilon(s,r)=((\cos(\pi/3),0,0,\sin(\pi/3));(1,(-1)^s,(-1)^r,(-1)^{s+r})/2),$$
where $r,s=0,1$, (the proof that $\dim\Fix<\epsilon>=2$ follows from Lemma
\ref{lem3}) or by an order two element, a plane reflection
$$\kappa(q,n,t)=(-1)^q((0,\cos(n\pi/3),\sin(n\pi/3),0);p^n((0,0,(-1)^t,1)/\sqrt{2})\,),$$
where $n=0,1,2$, $q,t=0,1$ and $p$ is the permutation $p(a,b,c,d)=(a,c,d,b)$.
The isotropy planes can be labelled as follows:
\begin{eqnarray*}\label{listex12}
P_1(s,r)=\Fix\Sigma_1(s,r),&\hbox{ where }\Sigma_1(s,r)=<\epsilon(s,r)>,\\
P_2(q,n,t)=\Fix\Sigma_2(q,n,t),&\hbox{ where }\Sigma_2(q,n,t)=<\kappa(q,n,t)>.
\end{eqnarray*}
Each of these plane contains exactly one copy of each of two types of
(not conjugate) symmetry axes, the isotropy groups of which are isomorphic
to $\D_3$. They are:
\begin{eqnarray*}\label{listex13}
L_1(s,r)&\hbox{ is the intersection of } P_1(s,r),\ P_2(0,0,s+1),
\ P_2(s+1,1,r)\hbox{ and }P_2(r,2,r+1),\\
L_2(s,r)&\hbox{ is the intersection of } P_1(s,r),\ P_2(1,0,s+1),
\ P_2(s,1,r)\hbox{ and }P_2(r+1,2,r+1).
\end{eqnarray*}
(Intersections of $P_1$ and $P_2$ or of two $P_2$'s follows from Lemma
\ref{lem2}.) Since $N(\Sigma_1)/\Sigma_1\cong\D_2$, the axes $L_1(s,r)$ and
$L_2(s,r)$ intersect orthogonally.
\\

\noindent Now, if
\begin{itemize}
\item[(i)] there exist two steady states $\xi_1\in L_1$ and $\xi_2\in L_2$;
\item[(ii)] $\xi_1$ is a saddle and $\xi_2$ is a sink in $P_1$,
moreover a saddle-sink heteroclinic orbit $\kappa_1$ exists between $\xi_1$ and $\xi_2$ in $P_1$;
\item[(iii)] $\xi_1$ is a sink and $\xi_2$ is a saddle in $P_2$,
moreover a saddle-sink heteroclinic orbit $\kappa_2$ exists between $\xi_1$ and $\xi_2$ in $P_1$,
\end{itemize}
then there exists a pseudo-simple robust heteroclinic cycle
$\xi_2\to\xi_1\to\xi_2$.
If in addition the unstable manifold of $\xi_2$ in $P_2$ is included
in the stable manifold of $\xi_1$: $(W^u(\xi_2)\cap P_1)\subset W^s(\xi_1)$,
then the system possesses a heteroclinic network comprised of two distinct
(i.e. not related by any symmetry) connections $\xi_2\to\xi_1$ and one
connection $\xi_1\to\xi_2$.
With minor modifications of Lemma 5 in \cite{pc13}, it can be proven
that the conditions (i)-(iii) above are satisfied for an open set of
$\Gamma$-equivariant vector fields. In the following subsection we
build an explicit system with these properties.

\subsection{Equivariant structure}\label{nform}

Here we derive a third order normal form commuting with the
action of the group $\Gamma=(\D_3\rl\Z_2;\mO\rl\V)$ introduced in subsection
\ref{d3v}. For convenience, we write the quaternion ${\bf q}=(q_1,q_2,q_3,q_4)$
as a pair of complex numbers,
\begin{equation}\label{cpair}
{\bf u}=(u_1,u_2),\hbox{ where }u_1=q_1+\ri q_2\hbox{ and }u_2=q_3+\ri q_4.
\end{equation}
The operation of multiplication (\ref{mqua}) takes the form
\begin{equation}\label{mcmpl}
{\bf h}{\bf u}=(h_1u_1-h_2\bar u_2,h_1u_2+h_2\bar u_1)
\end{equation}
and $\tilde{\bf u}=(\bar u_1,-u_2)$ is the conjugate of $\bf u$.
For quaternions presented as (\ref{cpair}) and points in $\R^4=\C^2$
as ${\bz}=(z_1,z_2)$, the action $(\bl;\br):\bz\to\bl\bz\br^{-1}$
of (some) elements of $\Gamma$ is
\begin{equation}\label{action}
\begin{array}{r|l}
(\bl;\br)&\bz\to\bl\bz\br^{-1}\\
\hline\\
((0,1,0,0);(1,1,0,0)/\sqrt{2}) &
(z_1,z_2)\to(\re^{\pi\ri/4}z_1,\re^{3\pi\ri/4}z_2)\\
\kappa(1,0,1)=((0,1,0,0);(0,0,1,-1)/\sqrt{2}) &
(z_1,z_2)\to(\re^{3\pi\ri/4}z_2,\re^{-3\pi\ri/4}z_1)\\
\kappa(0,0,1)=((0,1,0,0);(0,0,-1,1)/\sqrt{2}) &
(z_1,z_2)\to-(\re^{3\pi\ri/4}z_2,\re^{-3\pi\ri/4}z_1)\\
\epsilon(0,0)=((1,0,0,\sqrt{3})/2;(1,1,1,1)/2)&
(z_1,z_2)\to(\re^{-\pi\ri/4}z_1+\re^{-\pi\ri/4}z_2
-\ri\sqrt{3}(\re^{-\pi\ri/4}\bar z_1+\re^{-\pi\ri/4}\bar z_2;\\
&-\re^{\pi\ri/4}z_1+\re^{\pi\ri/4}z_2
+\ri\sqrt{3}(\re^{\pi\ri/4}\bar z_1+\re^{\pi\ri/4}\bar z_2).
\end{array}\end{equation}

Using the computer algebra program GAP, we derive from (\ref{action})
that the cubic equivariant terms are of the form
$$
\left(\begin{array}{c} \dot z_1 \\ \dot z_2 \end{array}\right)
= b(z_1\bar z_1+ z_2\bar z_2)\left(\begin{array}{c} z_1 \\ z_2 \end{array}
\right) + c\Theta_1+d\Theta_2+e\Theta_3
$$
where $b,c,d,e$ are real and the maps $\Theta_j$ are cubic expressions of
$z_1,\bar z_1, z_2,\bar z_2$:
$$
\Theta_1 = \left(\begin{array}{c} \sqrt{3}\bar z_1^2z_2 -\frac{\ri}{2}z_1^2\bar z_1+\ri z_1z_2\bar z_2-\frac{3\ri}{2}\bar z_1\bar z_2^2 \\ -\sqrt{3}z_1\bar z_2^2 -\frac{\ri}{2}z_2^2\bar z_2+\ri z_1\bar z_1 z_2-\frac{3 \ri}{2}\bar z_1^2\bar z_2 \end{array}\right)
$$
$$
\Theta_2 = \left(\begin{array}{c} z_2^3 - 3\bar z_1^2z_2 + 2\ri\sqrt{3}\bar z_1\bar z_2^2 \\ -z_1^3 +3z_1\bar z_2^2 + 2\ri\sqrt{3}\bar z_1^2\bar z_2 \end{array}\right)
$$
$$
\Theta_3 = \left(\begin{array}{c} -3\sqrt{3}z_1^2\bar z_1-\sqrt{3}\bar z_1\bar z_2^2+\ri z_2^3+3\ri\bar z_1^2z_2 \\ -3\sqrt{3}z_2^2\bar z_2-\sqrt{3}\bar z_1^2\bar z_2 - \ri z_1^3 - 3\ri z_1\bar z_2^2 \end{array}\right)
$$
In another arrangement of the cubic monomials, the equations read:
\begin{equation}\label{sys2}
\begin{array}{lcr}
\vspace*{2mm}
\dot z_1 &=& (\mu+b(z_1\bar z_1+ z_2\bar z_2))z_1+Az_1^2\bar z_1+\ri cz_1z_2\bar z_2+Bz_2^3+C\bar z_1^2z_2+D\bar z_1\bar z_2^2 \\
\dot z_2 &=& (\mu+b(z_1\bar z_1+ z_2\bar z_2))z_2 +Az_2^2\bar z_2+\ri cz_1\bar z_1z_2-Bz_1^3-Cz_1\bar z_2^2+D\bar z_1^2\bar z_2
\end{array}\end{equation}
where
\begin{equation}\label{coef}
A=-3\sqrt{3}e-\ri\frac{c}{2},\ B=d+\ri e,\ C=-3d+\sqrt{3}c+3\ri e,\
D=\sqrt{3}(-e+\ri(2d-\frac{\sqrt{3}}{2}c)).
\end{equation}
We set $\mu=1$ in order to make the origin unstable.

From (\ref{action}) we obtain that the invariant planes are:
\begin{align*}
P_2(1,0,1)&=\Fix(((0,1,0,0);(0,0,1,-1)/\sqrt{2}))=(\re^{3\pi\ri/4}w,w)\\
P_2(0,0,1)&=\Fix(((0,1,0,0);(0,0,-1,1)/\sqrt{2}))=(-\re^{3\pi\ri/4}w,w)\\
P_1(0,0)&=\Fix(((1,0,0,\sqrt{3})/2;(1,1,1,1)/2))=
(\sqrt{2}\re^{-3\pi\ri/4}w-\sqrt{3}\bar w,w)
\end{align*}
and the invariant axes are given by
$$L_1(0,0)=P_1(0,0)\cap P_2(0,0,1)=(-\re^{3\pi\ri/4}r_1\re^{\alpha_1},r_1\re^{\alpha_1}),
\hbox{ where }\re^{2\alpha_1}={-1+\sqrt{2}\ri\over\sqrt{3}}\re^{\pi\ri/4},$$
and
$$L_2(0,0)=P_1(0,0)\cap P_2(1,0,1)=(\re^{3\pi\ri/4}r_2\re^{\alpha_2},r_2\re^{\alpha_2}),
\hbox{ where } \re^{2\alpha_2}={1+\sqrt{2}\ri\over\sqrt{3}}\re^{\pi\ri/4}.$$
The system (\ref{sys2}) restricted onto the axes $L_1$ is
$$\dot r_1=
r_1(\mu+r_1^2(2b+2\sqrt{2}c-{8\sqrt{2}\over\sqrt{3}}d-{4\over\sqrt{3}}e)),$$
and onto $L_2$
$$\dot r_2=
r_2(\mu+r_2^2(2b-2\sqrt{2}c+{8\sqrt{2}\over\sqrt{3}}d-{4\over\sqrt{3}}e)).$$
Therefore, the conditions for existence of the steady states $\xi_1\in L_1$ and $\xi_2\in L_2$ are
\begin{equation}\label{cond1}
2b+2\sqrt{2}c-{8\sqrt{2}\over\sqrt{3}}d-{4\over\sqrt{3}}e<0\hbox{ and }
2b-2\sqrt{2}c+{8\sqrt{2}\over\sqrt{3}}d-{4\over\sqrt{3}}e<0,
\end{equation}
respectively. Whenever the steady states exist, the eigenvalues of $df(\xi_j)$
in the radial directions are $-2\mu$.

By considering the restriction of the system into the subspace
$\Fix(-((0,1,0,0);(0,0,1,-1)/\sqrt{2}))$ we obtain that the eigenvalues
of $df(\xi_1)$ and $df(\xi_2)$, associated with two-dimensional eigenspaces,
are
\begin{equation}\label{cond2}
\lambda_1^{\rm mlt}={2\sqrt{2}r_1^2\over3}(-9c+10\sqrt{3}d-2\sqrt{6}e)
~~\hbox{ and }~~
\lambda_2^{\rm mlt}={2\sqrt{2}r_2^2\over3}(9c-10\sqrt{3}d-2\sqrt{6}e).
\end{equation}
By considering the restriction of the system into the subspace fixed by
$((1,0,0,\sqrt{3})/2;(1,1,1,1)/2)$ we derive the expressions for the remaining
eigenvalues:
\begin{equation}\label{cond3}
\lambda_1^{\rm sgl}=r_1^2(-{32\over\sqrt{3}}e+{8\sqrt{2}\over\sqrt{3}}d)
~~\hbox{ and }~~
\lambda_2^{\rm sgl}=r_2^2(-{32\over\sqrt{3}}e-{8\sqrt{2}\over\sqrt{3}}d).
\end{equation}
Necessary conditions for the existence of heteroclinic cycle $\xi_1\to\xi_2\to\xi_1$
discussed in subsection \ref{d3v} are that
\begin{equation}\label{cexist}
\lambda_1^{\rm sgl}>0,~~~\lambda_2^{\rm sgl}<0,~~~
\lambda_1^{\rm mlt}<0~\hbox{ and }\lambda_2^{\rm mlt}>0.
\end{equation}

\subsection{The numerical simulations}

The general third order $\Gamma$-equivariant system is given by
(\ref{sys2})-(\ref{coef}). From (\ref{cond1})-(\ref{cexist}), the
necessary conditions for existence of the heteroclinic cycle
expressed in terms of the coefficients $b$, $c$, $d$ and $e$ are:
$$-d<2\sqrt{2}e<d,\
c>{1\over9}\max((10\sqrt{3}d+2\sqrt{6}e),(10\sqrt{3}d-2\sqrt{6}e)),$$
$$b<{1\over3}\min((3\sqrt{2}c-4\sqrt{6}d+2\sqrt{3}e),(-3\sqrt{2}c+4\sqrt{6}d+2\sqrt{3}e)).$$
We set
$$d=1,\ e={1\over2\sqrt{2}}(d-h_1),\ c={1\over9}(10\sqrt{3}d+2\sqrt{6}e)+h_2,\
b={1\over3}(3\sqrt{2}c-4\sqrt{6}d+2\sqrt{3}e)-1,$$
$$\hbox{where }~-1<h_1<1~\hbox{ and }~h_2~\hbox{ is positive.}$$
In Section \ref{sec4} we have proved that in a $\Gamma$-equivariant system,
for the considered group $\Gamma\cong\D_3$, an asymptotically stable periodic orbit
bifurcates from a pseudo-simple heteroclinic cycle when the double expanding
eigenvalue vanishes. It can be easily shown that the expressions for local and
global maps near a heteroclinic cycle in a $(\D_3\rl\Z_2;\mO\rl\V)$-equivariant
system are identical to the ones given in the proof of Theorem \ref{thperorb},
with $\gamma\in(\D_3\rl\Z_2;\mO\rl\V)$ in (\ref{glmap1n}).
Therefore, the proof holds true for $\Gamma=(\D_3\rl\Z_2;\mO\rl\V)$ that we
consider here. The expanding double eigenvalue is
$\lambda_2^{\rm mlt}=2\sqrt{2}h_2r_1^2/3$, hence we consider small values of $h_2$.

\begin{figure}[p]
\hspace*{-2cm}\includegraphics[width=10cm]{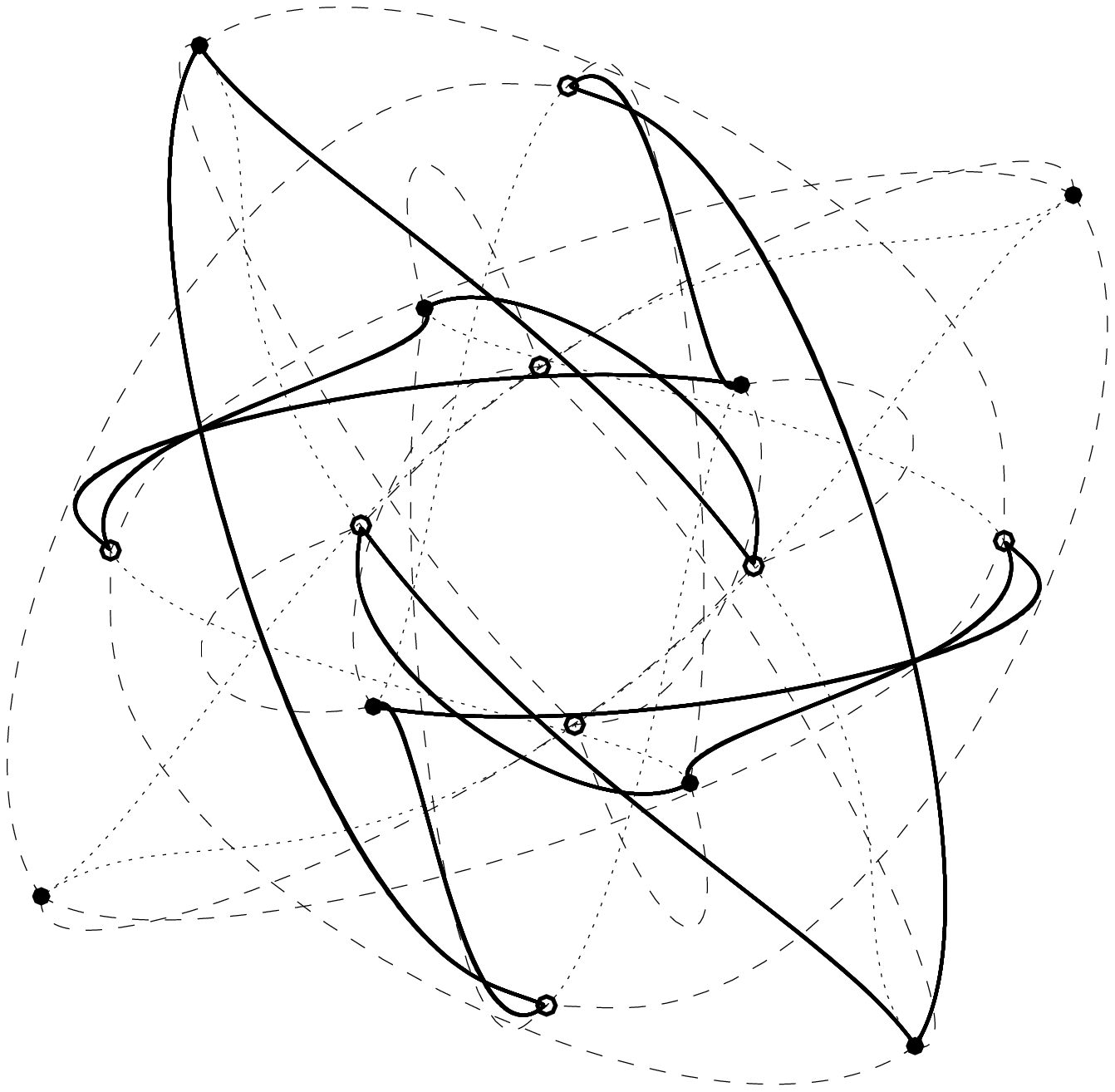}

\vspace*{-15mm}
\hspace*{2cm}{\large (a)}

\vspace*{-39mm}
\hspace*{6cm}\includegraphics[width=10cm]{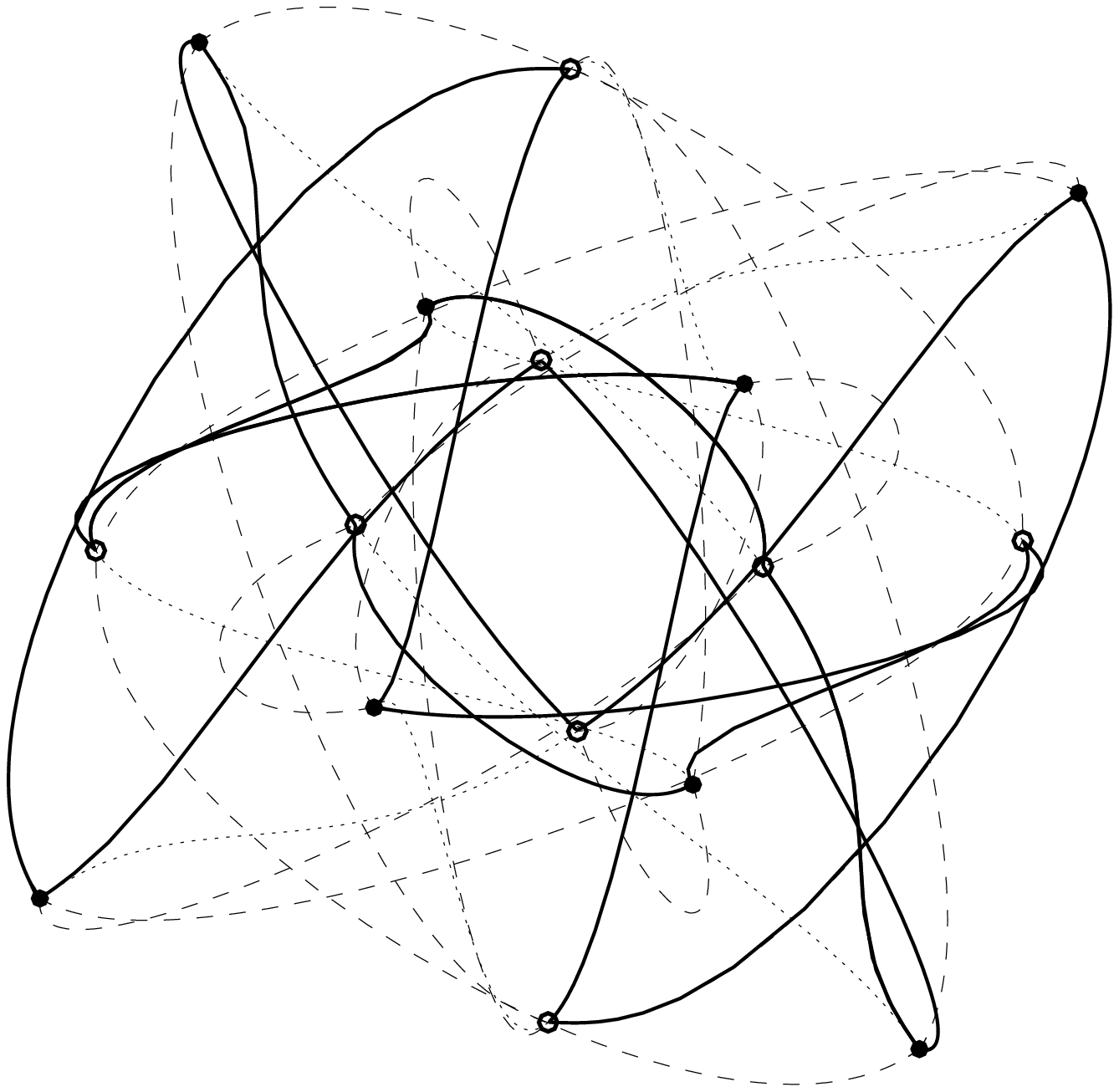}

\vspace*{-2cm}
\hspace*{10cm}{\large (b)}

\vspace*{-29mm}
\hspace*{-2cm}\includegraphics[width=10cm]{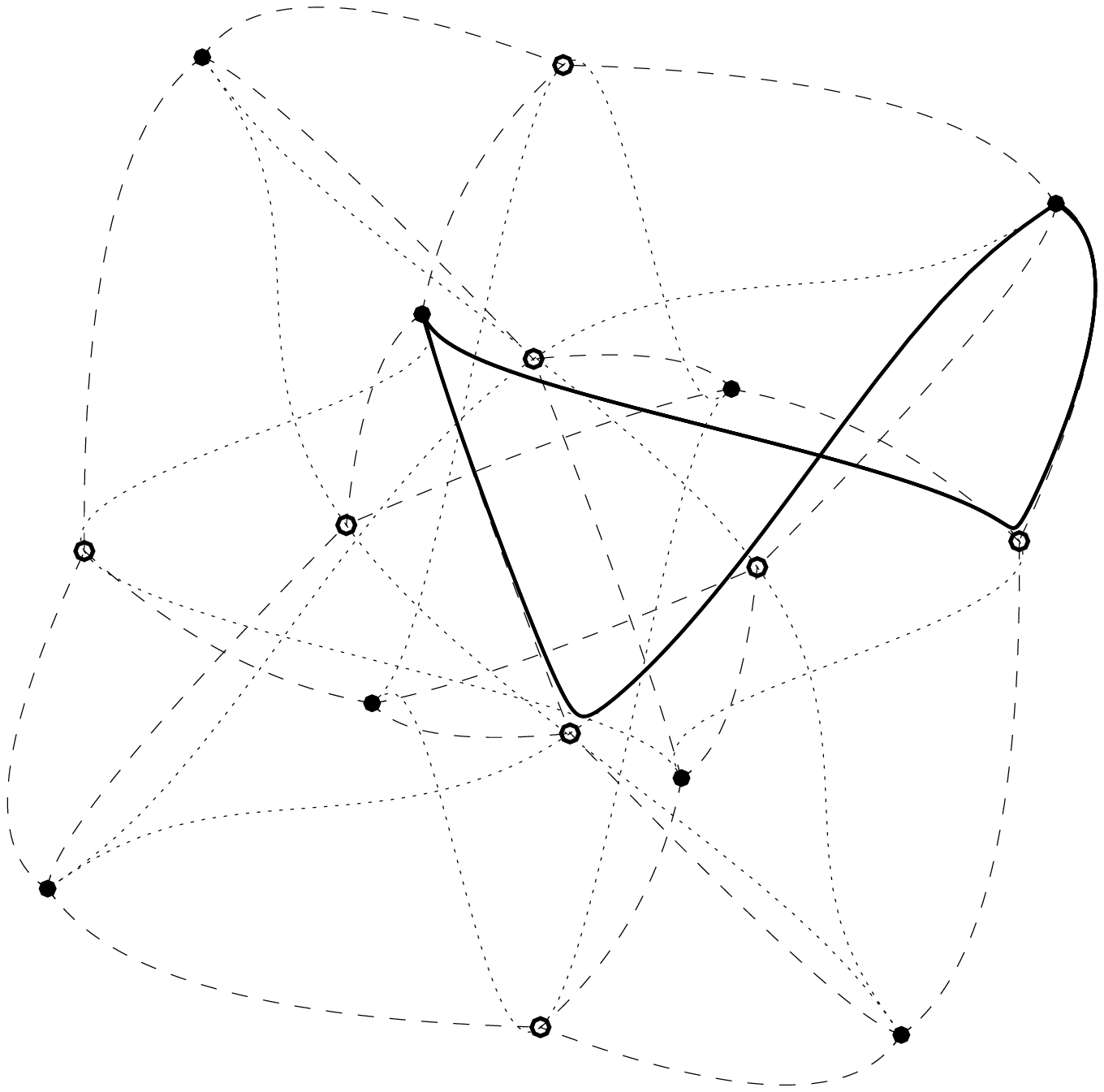}

\vspace*{-16mm}
\hspace*{2cm}{\large (c)}

\vspace*{4mm}
\noindent
\caption{Projection of the periodic orbits (solid line) and of the
heteroclinic connections (dashed and dotted lines) into the plane
$<{\bf v}_1,{\bf v}_2>$, where ${\bf v}_1=(1,2,1,1.8)$ and
${\bf v}_2=(2,-1,1.8,-1)$, for $h_2=0.001$ and $h_1=0.7$ (a), $h_1=0.8$ (b)
and $h_1=0.92$ (c). Steady states $\xi_1$ are denoted by hollow circles and
$\xi_2$ by filled ones.}
\label{fig4}\end{figure}

\begin{figure}[p]

\vspace*{-14mm}
\hspace*{2mm}\includegraphics[height=30mm,width=16cm]{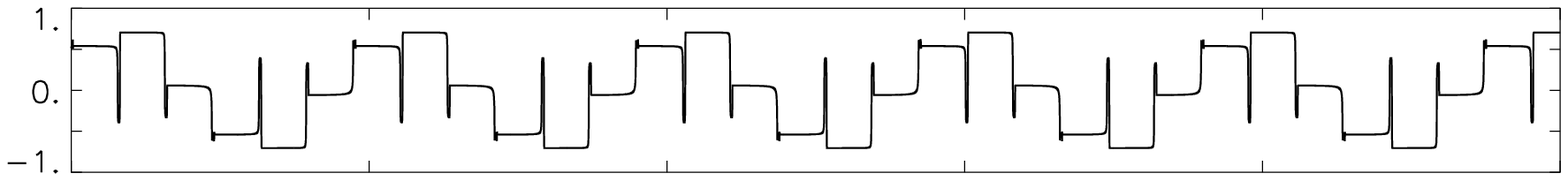}

\vspace*{-15mm}
\hspace*{2mm}\includegraphics[height=30mm,width=16cm]{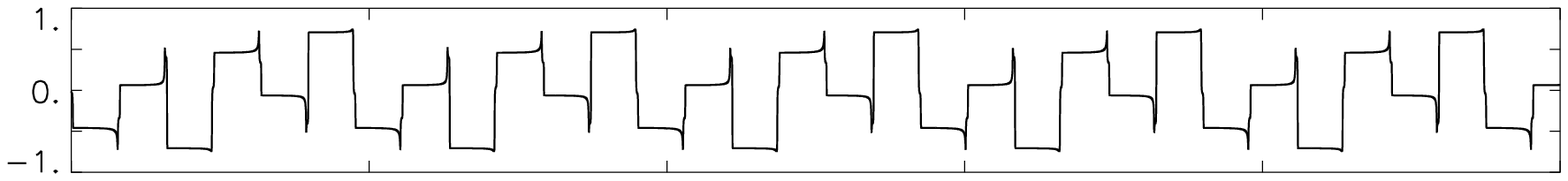}

\vspace*{-15mm}
\hspace*{2mm}\includegraphics[height=30mm,width=16cm]{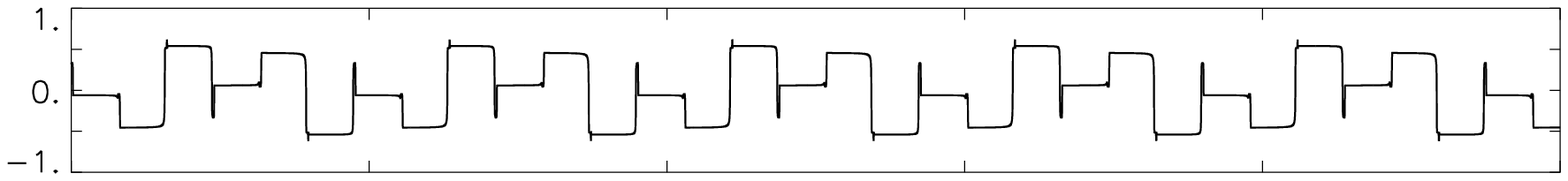}

\vspace*{-15mm}
\hspace*{2mm}\includegraphics[height=30mm,width=16cm]{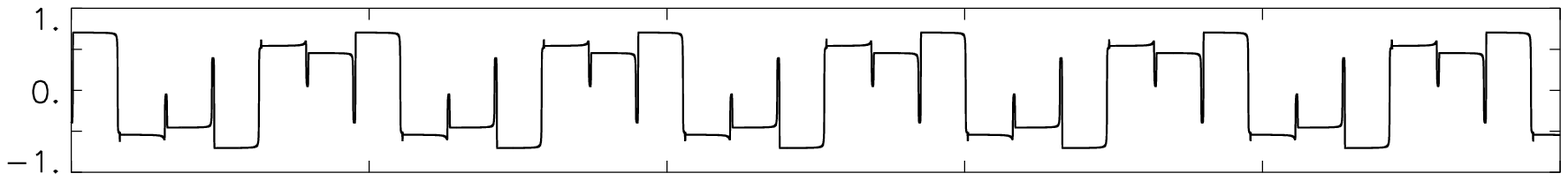}

\vspace*{-11mm}
\hspace*{15mm}{10000\hspace*{125mm}20000}

\vspace*{-3mm}
\hspace*{122mm}{\large $t$}

\vspace*{-64mm}Re($z_1$)

\vspace*{10mm}Im($z_1$)

\vspace*{10mm}Re($z_2$)

\vspace*{10mm}Im($z_2$)

\vspace*{10mm}
\hspace*{53mm}{\large (a)}

\hspace*{2mm}\includegraphics[height=30mm,width=16cm]{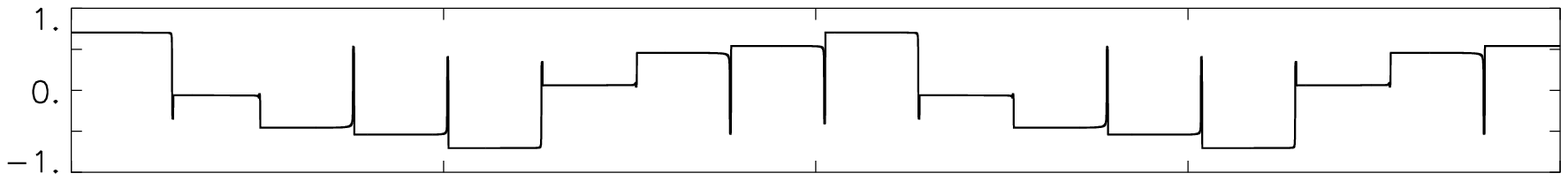}

\vspace*{-15mm}
\hspace*{2mm}\includegraphics[height=30mm,width=16cm]{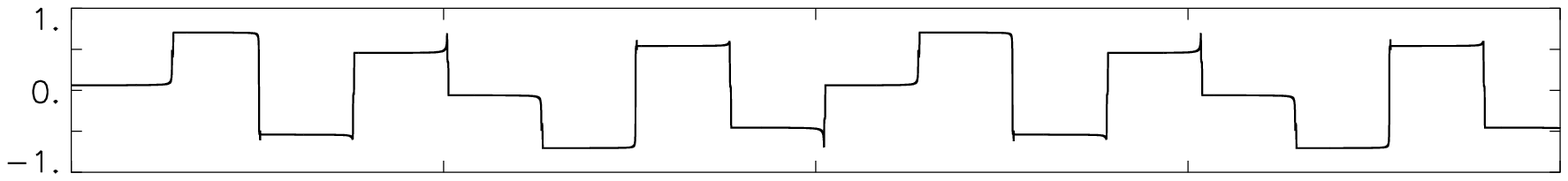}

\vspace*{-15mm}
\hspace*{2mm}\includegraphics[height=30mm,width=16cm]{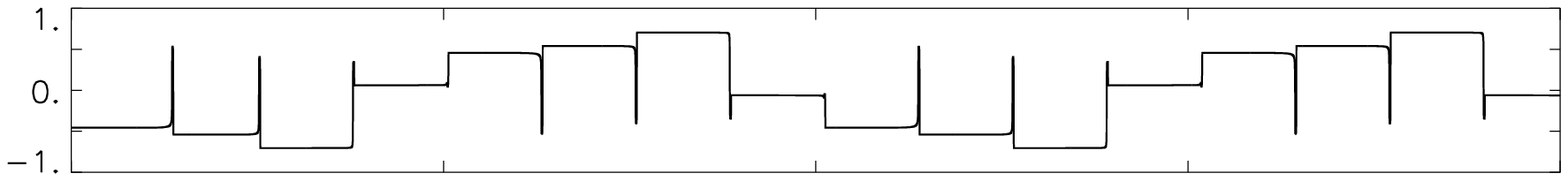}

\vspace*{-15mm}
\hspace*{2mm}\includegraphics[height=30mm,width=16cm]{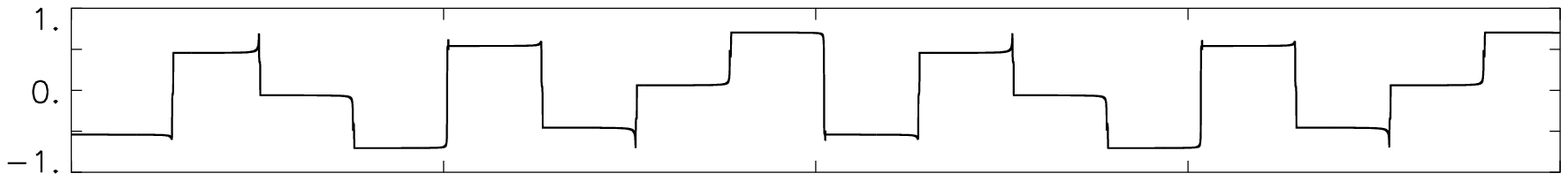}

\vspace*{-11mm}
\hspace*{21mm}{0\hspace*{128mm}20000}

\vspace*{-3mm}
\hspace*{122mm}{\large $t$}

\vspace*{-64mm}Re($z_1$)

\vspace*{10mm}Im($z_1$)

\vspace*{10mm}Re($z_2$)

\vspace*{10mm}Im($z_2$)

\vspace*{10mm}
\hspace*{53mm}{\large (b)}

\hspace*{2mm}\includegraphics[height=30mm,width=16cm]{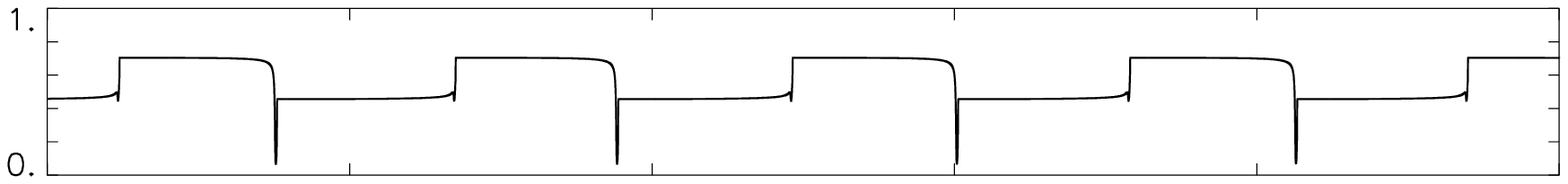}

\vspace*{-15mm}
\hspace*{2mm}\includegraphics[height=30mm,width=16cm]{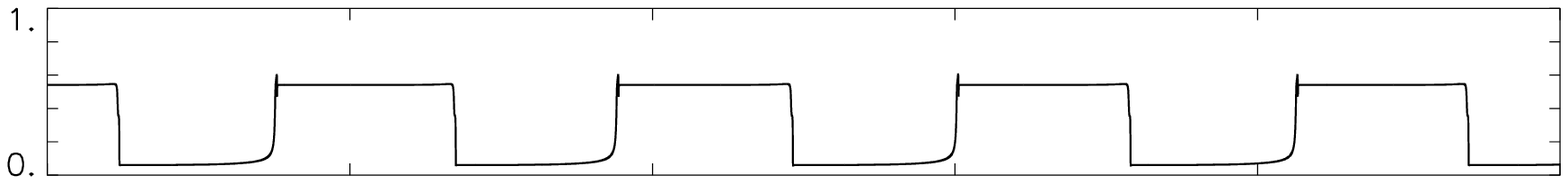}

\vspace*{-15mm}
\hspace*{2mm}\includegraphics[height=30mm,width=16cm]{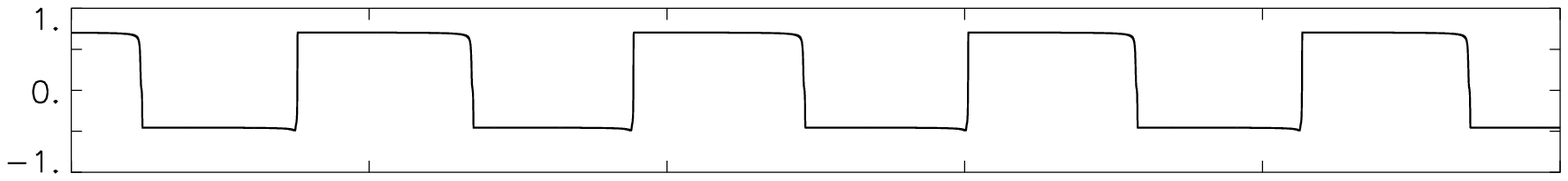}

\vspace*{-15mm}
\hspace*{2mm}\includegraphics[height=30mm,width=16cm]{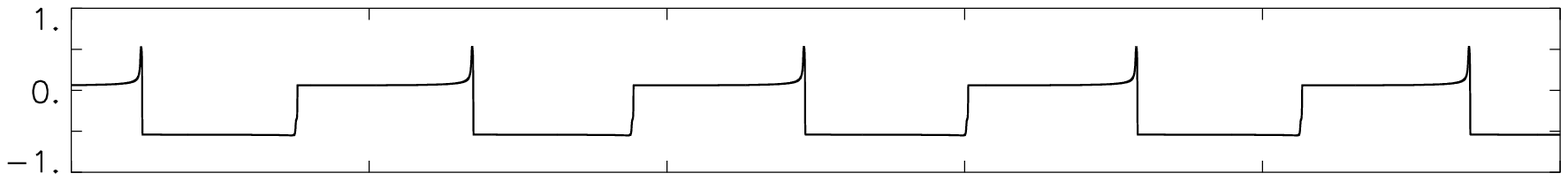}

\vspace*{-11mm}
\hspace*{15mm}{10000\hspace*{125mm}20000}

\vspace*{-3mm}
\hspace*{122mm}{\large $t$}

\vspace*{-64mm}Re($z_1$)

\vspace*{10mm}Im($z_1$)

\vspace*{10mm}Re($z_2$)

\vspace*{10mm}Im($z_2$)

\vspace*{10mm}
\hspace*{53mm}{\large (c)}

\noindent
\caption{Temporal behaviour of $z_1$ and $z_2$
for $h_2=0.001$ and $h_1=0.7$ (a), $h_1=0.8$ (b) and $h_1=0.92$ (c).}
\label{fig5}\end{figure}

The computations were done for several values of $h_1$, $0.5<h_1<1$.
For small enough $h_2$ the attractor is a periodic orbit near
the group orbit of heteroclinic connections $\kappa_1$ and $\kappa_2$,
$$X=(\cup_{\gamma\in\Gamma}\gamma\kappa_1)\bigcup
(\cup_{\gamma\in\Gamma}\gamma\kappa_2).$$
Three instances of these orbits are shown in Figs.~\ref{fig4} and \ref{fig5}.
Note, that by varying $h_1$ we get periodic orbits with different $\gamma$'s
in (\ref{glmap1n}) (see remark \ref{varorb}).

\begin{figure}[p]
\hspace*{-2cm}\includegraphics[width=10cm]{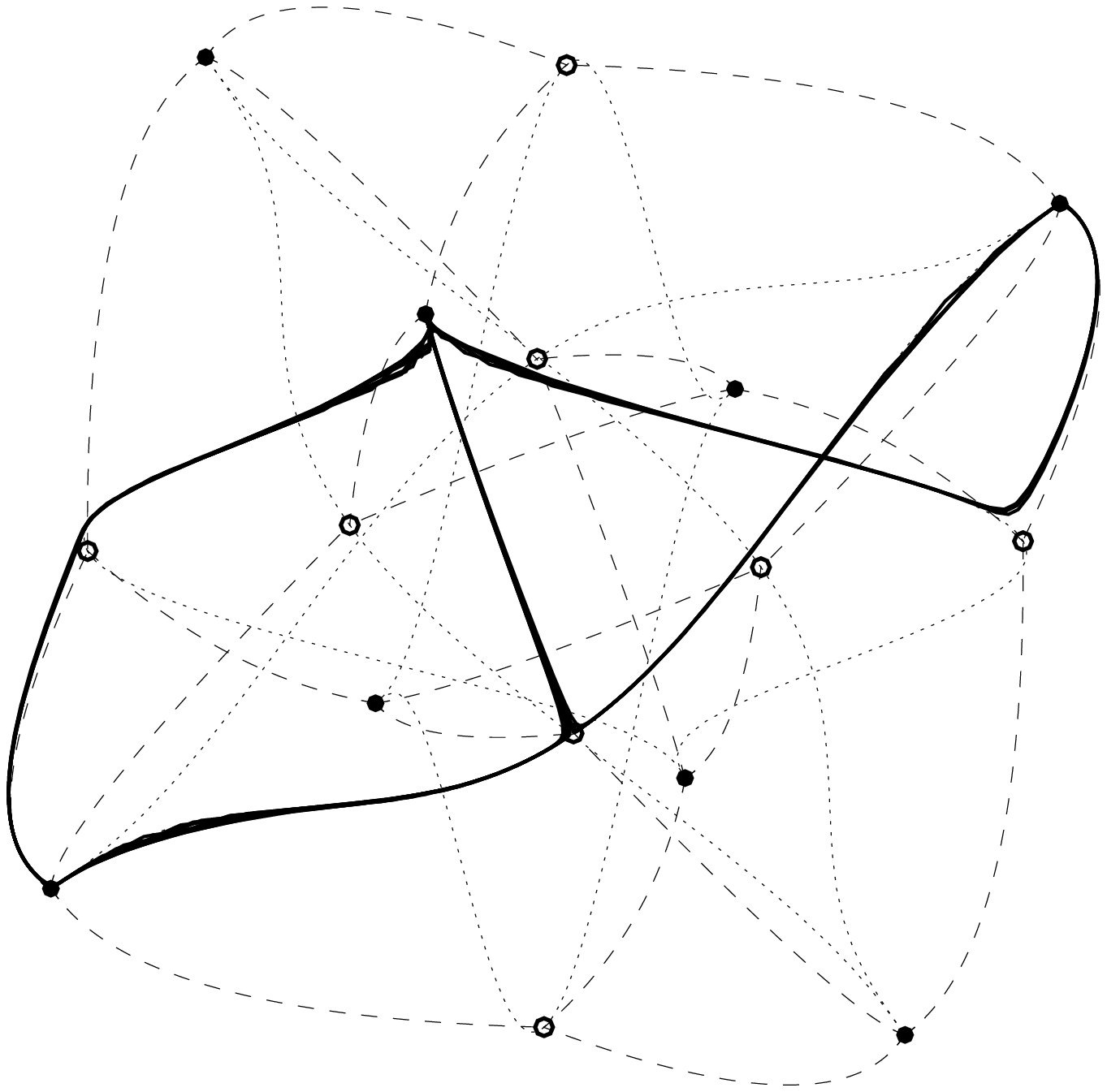}

\vspace*{-15mm}
\hspace*{2cm}{\large (a)}

\vspace*{-39mm}
\hspace*{6cm}\includegraphics[width=10cm]{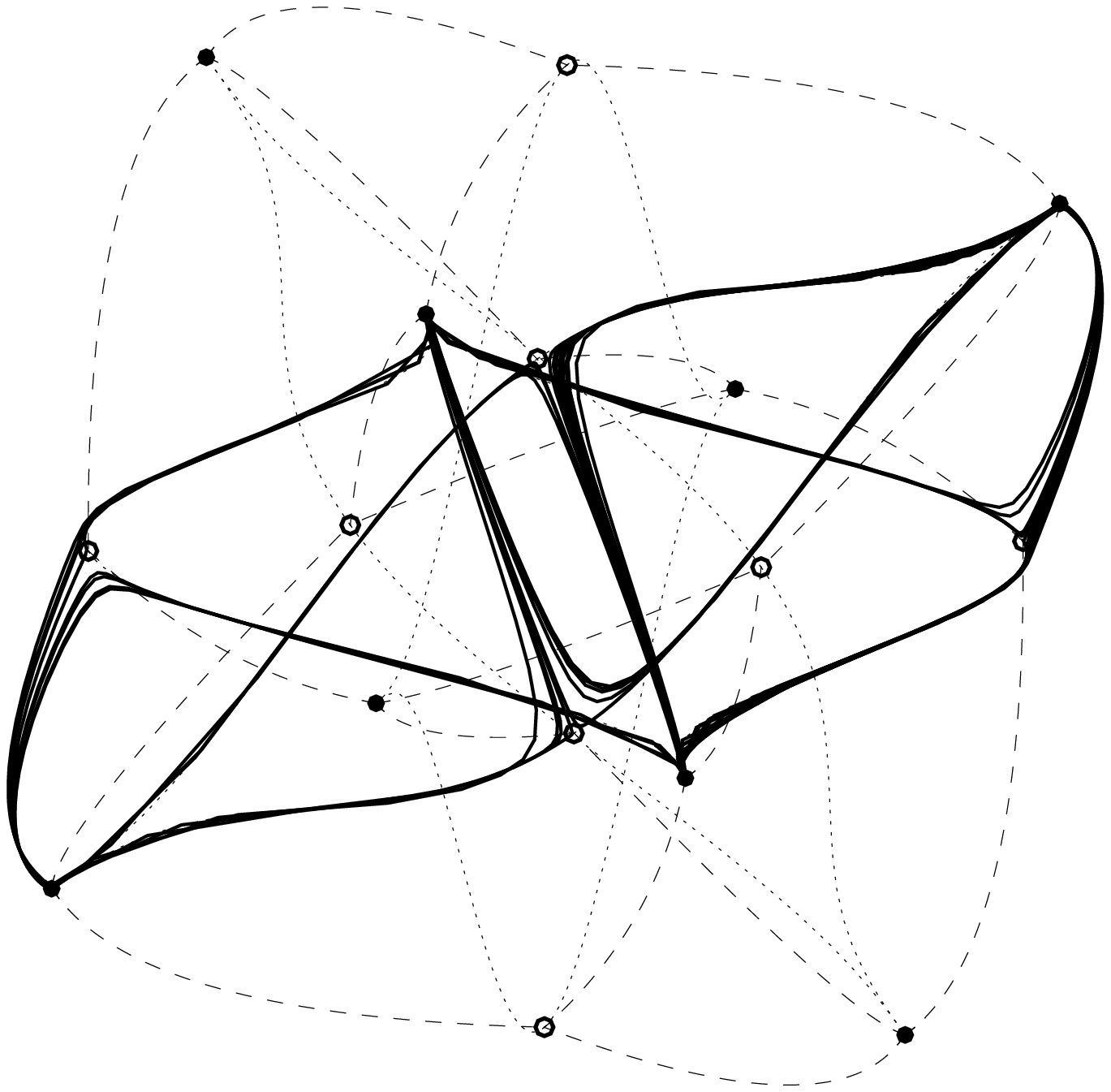}

\vspace*{-16mm}
\hspace*{10cm}{\large (b)}

\vspace*{-29mm}
\hspace*{-2cm}\includegraphics[width=10cm]{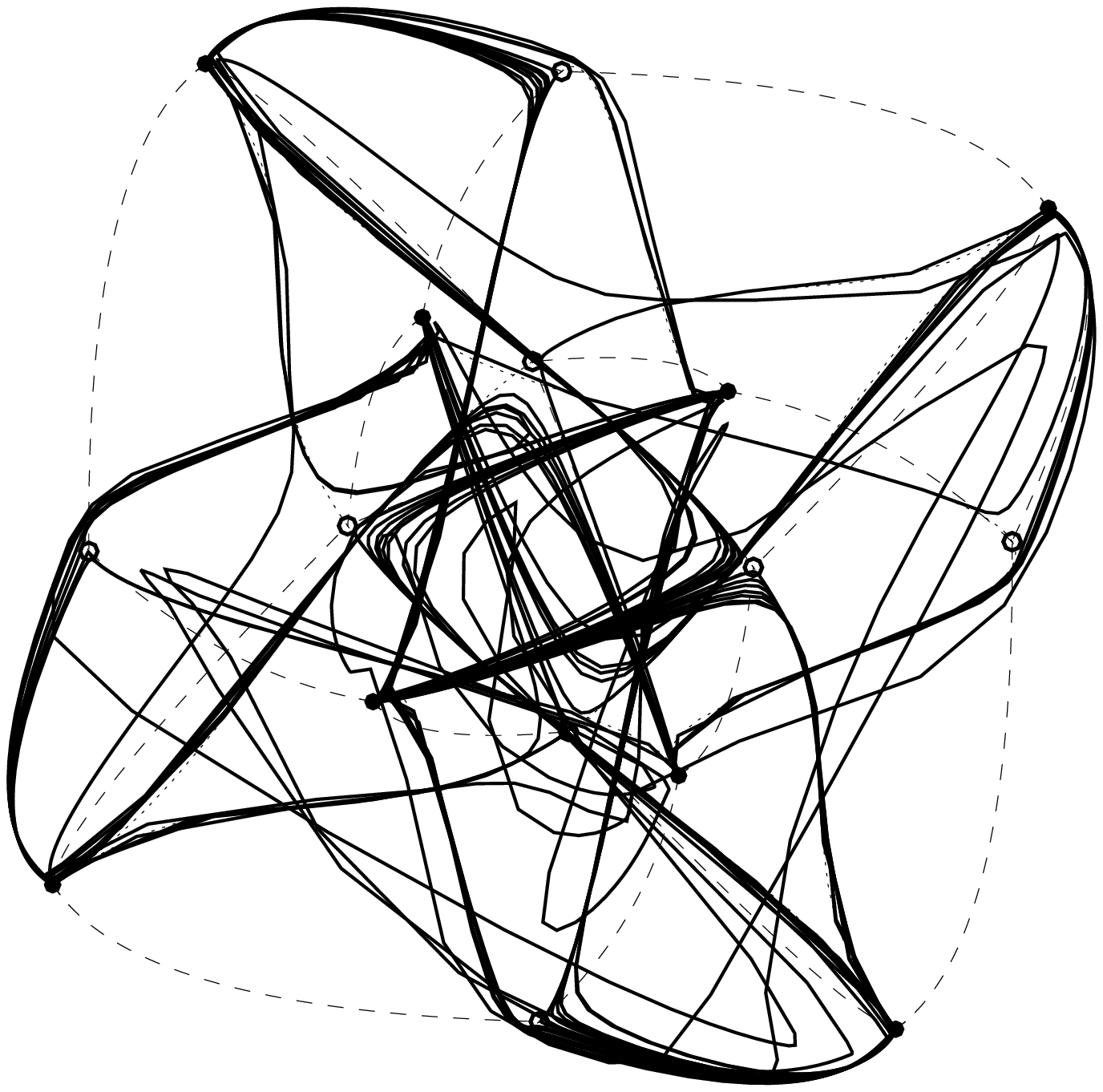}

\vspace*{-16mm}
\hspace*{2cm}{\large (c)}

\vspace*{4mm}
\noindent
\caption{Projection of the chaotic trajectories (solid line) and of the
heteroclinic connections (dashed and dotted lines) into the plane
$<{\bf v}_1,{\bf v}_2>$, where ${\bf v}_1=(1,2,1,1.8)$ and
${\bf v}_2=(2,-1,1.8,-1)$, for $h_1=0.92$ and $h_2=0.0015$ (a), $h_2=0.002$ (b)
and $h_2=0.0028$ (c). Steady states $\xi_1$ are denoted by hollow circles and
$\xi_2$ by filled ones.}
\label{fig6}\end{figure}

\begin{figure}[p]

\vspace*{-14mm}
\hspace*{2mm}\includegraphics[height=30mm,width=16cm]{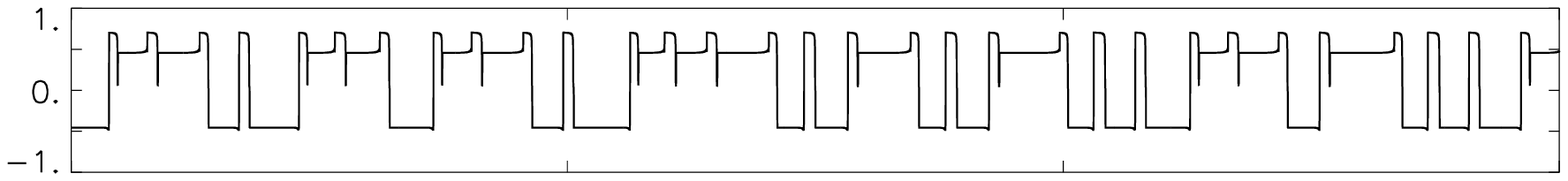}

\vspace*{-15mm}
\hspace*{2mm}\includegraphics[height=30mm,width=16cm]{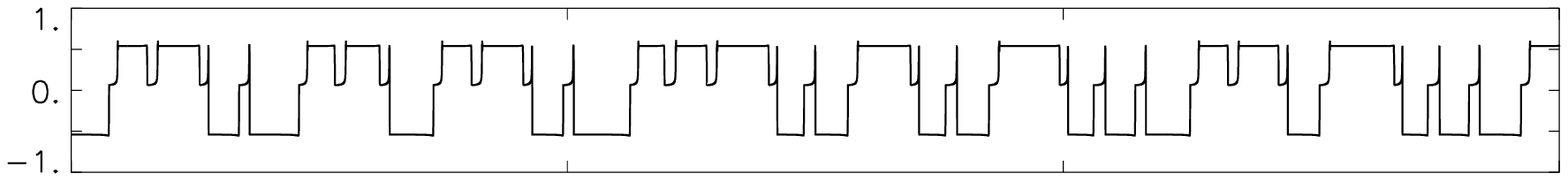}

\vspace*{-15mm}
\hspace*{2mm}\includegraphics[height=30mm,width=16cm]{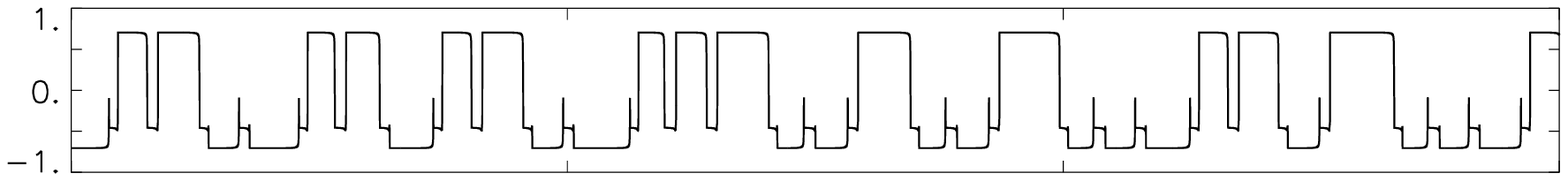}

\vspace*{-15mm}
\hspace*{2mm}\includegraphics[height=30mm,width=16cm]{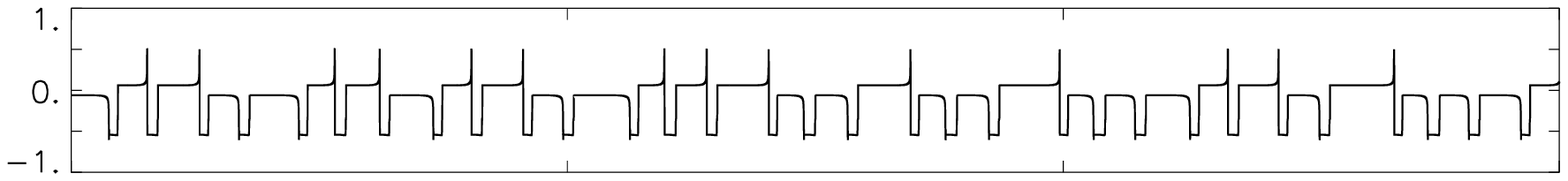}

\vspace*{-11mm}
\hspace*{15mm}{20000\hspace*{125mm}80000}

\vspace*{-3mm}
\hspace*{122mm}{\large $t$}

\vspace*{-64mm}Re($z_1$)

\vspace*{10mm}Im($z_1$)

\vspace*{10mm}Re($z_2$)

\vspace*{10mm}Im($z_2$)

\vspace*{10mm}
\hspace*{53mm}{\large (a)}

\hspace*{2mm}\includegraphics[height=30mm,width=16cm]{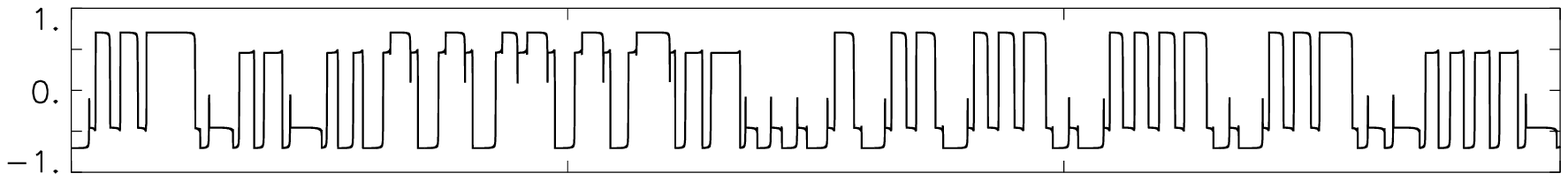}

\vspace*{-15mm}
\hspace*{2mm}\includegraphics[height=30mm,width=16cm]{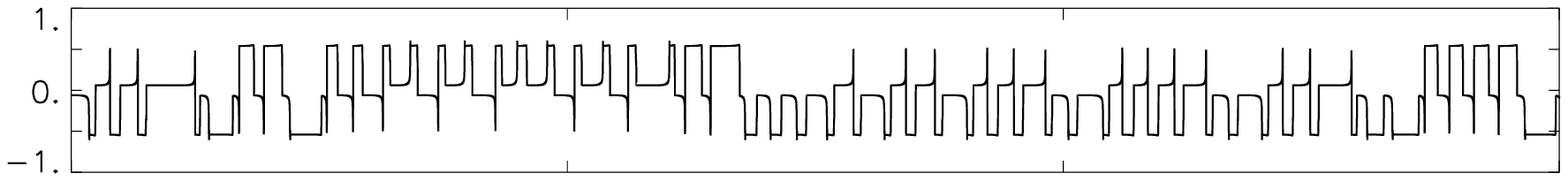}

\vspace*{-15mm}
\hspace*{2mm}\includegraphics[height=30mm,width=16cm]{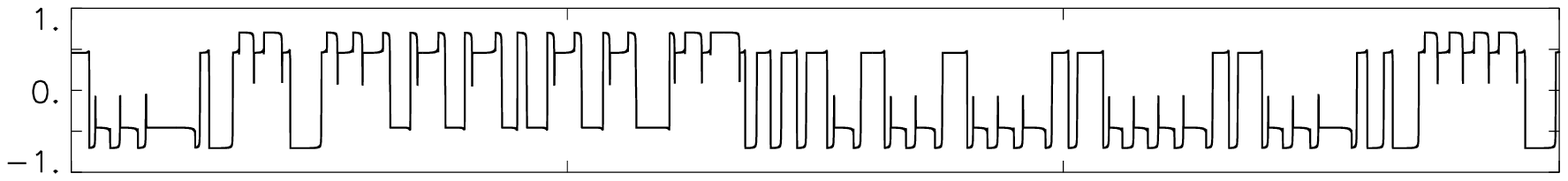}

\vspace*{-15mm}
\hspace*{2mm}\includegraphics[height=30mm,width=16cm]{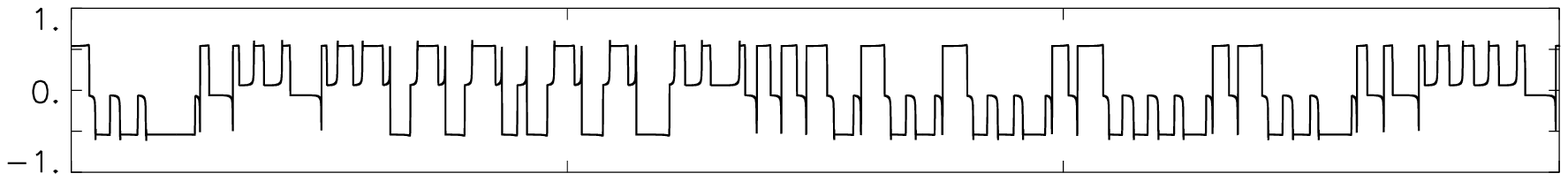}

\vspace*{-11mm}
\hspace*{15mm}{20000\hspace*{125mm}80000}

\vspace*{-3mm}
\hspace*{122mm}{\large $t$}

\vspace*{-64mm}Re($z_1$)

\vspace*{10mm}Im($z_1$)

\vspace*{10mm}Re($z_2$)

\vspace*{10mm}Im($z_2$)

\vspace*{10mm}
\hspace*{53mm}{\large (b)}

\hspace*{2mm}\includegraphics[height=30mm,width=16cm]{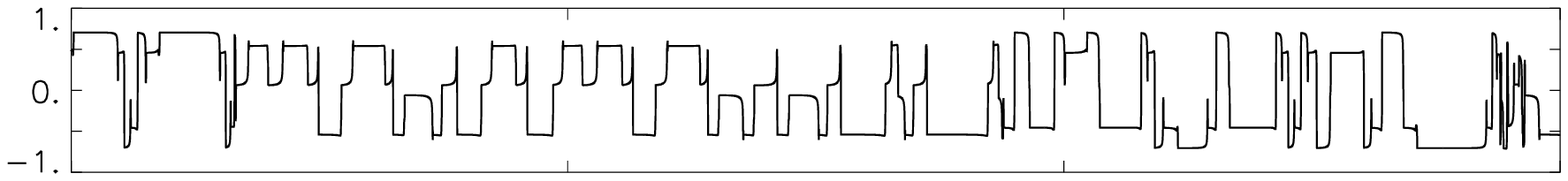}

\vspace*{-15mm}
\hspace*{2mm}\includegraphics[height=30mm,width=16cm]{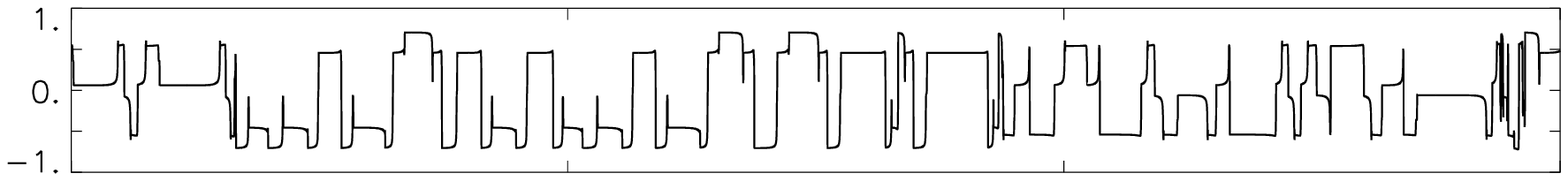}

\vspace*{-15mm}
\hspace*{2mm}\includegraphics[height=30mm,width=16cm]{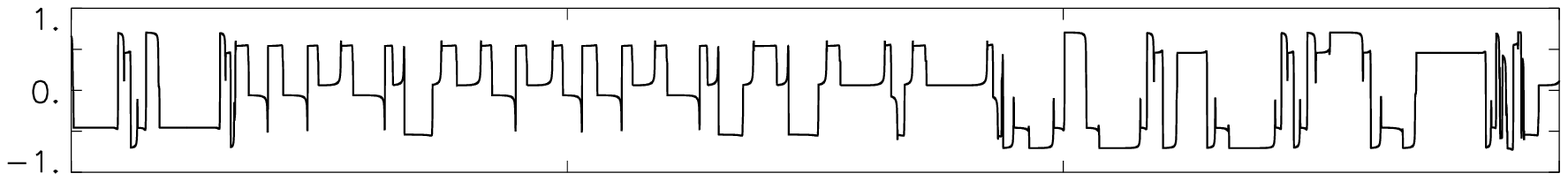}

\vspace*{-15mm}
\hspace*{2mm}\includegraphics[height=30mm,width=16cm]{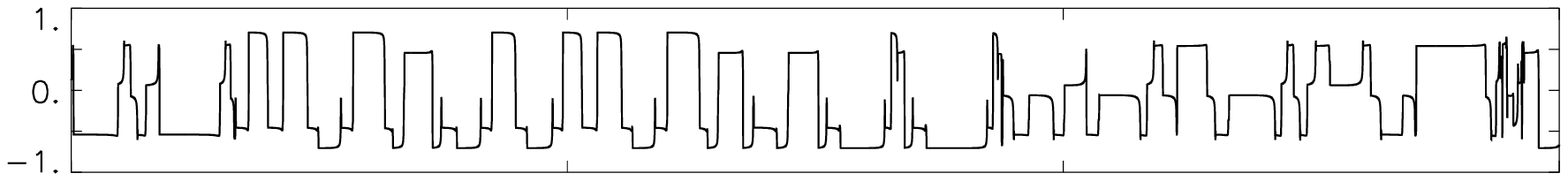}

\vspace*{-11mm}
\hspace*{15mm}{50000\hspace*{125mm}80000}

\vspace*{-3mm}
\hspace*{122mm}{\large $t$}

\vspace*{-64mm}Re($z_1$)

\vspace*{10mm}Im($z_1$)

\vspace*{10mm}Re($z_2$)

\vspace*{10mm}Im($z_2$)

\vspace*{10mm}
\hspace*{53mm}{\large (c)}

\noindent
\caption{Temporal behaviour of $z_1$ and $z_2$
for $h_1=0.92$ and $h_2=0.0015$ (a), $h_2=0.002$ (b) and $h_2=0.0028$ (c).}
\label{fig7}\end{figure}

With the increase of $h_2$ the behaviour ceases to be time-periodic,
however the chaotic trajectory stays close to the heteroclinic cycle. The
maximal distance of the points of the trajectory from the heteroclinic
cycle increases with $h_2$ (see Figs.~\ref{fig6} and \ref{fig7}).

\section{Conclusion}\label{conclusion}

Heteroclinic cycles in $\R^4$ have been extensively studied in the last
three decades, however only recently \cite{pc13} existence of pseudo-simple
heteroclinic cycles has been noticed. The present paper is the first one
where properties of these cycles are investigated.
Our main result is that pseudo-simple heteroclinic cycles in $\R^4$ are
{\em completely unstable} under generic conditions when the symmetry group
$\Gamma$ of the system is contained in SO(4). We further analysed the asymptotic
behavior near a pseudo-simple cycle with a specific group $\Gamma$ such that
the heteroclinic cycle connects precisely two
different equilibria. We proved that: (i) when the unstable double eigenvalue
at one equilibrium is close to 0, a nearby stable periodic orbit can exist
under generic conditions, (ii) when $\Gamma$ is extended to a group
$\tilde\Gamma$ which is not a subgroup of $SO(4)$, then the pseudo-simple
heteroclinic cycles are {\em fragmentarily asymptotically stable}. These
properties have been illustrated numerically. Finally, in the last part of
this study we have considered a more complex subgroup of SO(4), for which the
proof of existence of periodic orbits near the heteroclinic cycles is still
valid but which possess a richer structure. Numerical simulations have shown
that periodic orbits can follow different connections along the group of
heteroclinic cycle and that non-periodic attractors in
the vicinity of the cycle can also exist.

We expect that the results of Section \ref{sec4} can
be generalized to other subgroups $\Gamma\subset$\,SO(4), at least when
the unstable two-dimensional manifolds are invariant by the action of $\D_3$.
We also intend to look at the case when the two-dimensional unstable manifolds
are invariant under the action of $\D_k$ with $k>3$, the considered
example being readily modified to $\Gamma\cong\D_k$.
Arguments of Theorem \ref{thsect5} are likely to hold true for
other subgroups of O(4), which are not in SO(4), admitting
pseudo-simple heteroclinic cycles.
In light of our numerical observations in Section 6, we think it will be of
interest to investigate further the transition to complex dynamics near
a pseudo-simple cycle with symmetry group in SO(4), for groups studied in
this paper, of for different ones.

The definition of pseudo-simple heteroclinic cycles can be generalized
to $\R^n$ with $n>4$ by requiring the unstable eigenvalue at one of
the equilibria to have dimension of the associated eigenspace to be larger than
one. Evidently, such cycles are not asymptotically stable. The question whether
they can be f.a.s. for $\Gamma\subset$\,SO($n$) is yet another open problem.

\end{document}